\documentclass[letterpaper, 10 pt, conference]{ieeeconf}  % Comment this line out if you need a4paper

\IEEEoverridecommandlockouts                              % This command is only needed if 
\usepackage[fleqn]{amsmath}
\usepackage{amsfonts}
\usepackage{mathtools}
\usepackage{commath}
\usepackage{amssymb}
\usepackage{multicol}
\usepackage{breqn}
\usepackage{cite}
\usepackage{tabu}
\usepackage{graphicx}
\usepackage{subcaption}
\usepackage{url}
\usepackage{color}
\usepackage[bottom]{footmisc}
\graphicspath{{figures/}}

\newtheorem{thm}{Theorem}
\newtheorem{assumption}{Assumption}
\newtheorem{lem}{Lemma}

\title{\LARGE \bf
Adaptive Hessian Estimation Based Extremum Localization
}

\author{Huseyin Demircioglu, Iman Fadakar, Baris Fidan% <-this % stops a space
\thanks{*This work is supported by the Canadian NSERC Discovery Grant 116806.}% <-this % stops a space
\thanks{Huseyin Demircioglu, Iman Fadakar and Baris Fidan are with The Department of Mechanical and Mechatronics Engineering , University of Waterloo, Waterloo, ON, N2L 3G1, Canada hdemircioglu@uwaterloo.ca; ifadakar@uwaterloo.ca; fidan@uwaterloo.ca}%
}

\begin{document}

\maketitle
\thispagestyle{empty}
\pagestyle{empty}

%%%%%%%%%%%%%%%%%%%%%%%%%%%%%%%%%%%%%%%%%%%%%%%%%%%%%%%%%%%%%%%%%%%%%%%%%%%%%%%%
\begin{abstract}

In this paper we study continuous time adaptive extremum localization of an arbitrary quadratic function $F(\cdot)$ based on Hessian estimation, using measured the signal intensity by a sensory agent. The function $F(\cdot)$ represents a signal field as a result of a source located at the maximum point of $F(\cdot)$ and is decreasing as moving away from the source location. Stability of the proposed adaptive estimation and localization scheme is analyzed and the Hessian parameter and location estimates are shown to asymptotically converge to the true values. Moreover, the stability and convergence properties of algorithm are shown to be robust to drift in the extremum location. Simulation test results are displayed to verify the established properties of the proposed scheme as well as robustness to signal measurement noise.  
\end{abstract}

%%%%%%%%%%%%%%%%%%%%%%%%%%%%%%%%%%%%%%%%%%%%%%%%%%%%%%%%%%%%%%%%%%%%%%%%%%%%%%%%
\section{INTRODUCTION}

Years by years, source localization  yields some promiser applications, hence, it has been studied broadly such as \cite{mao2009localization, umay2017adaptive, sayed2005network,bishop2010optimality}. The generic task in these problems is that one or more sensory agents locate the source of a signal field with the help of measurement obtained from sensors mounted on these agents.  In order to localize the source, different kinds of measurements are utilized depending on the on the setting and constraints of the particular localization task. Generally, localization is accomplished using the information of the relative position of a single agent or multi-agents to a source such as bearing / angle of arrival (AOA) \cite{niculescu2004positioning,klukas1998line}, time difference of arrival (TDOA) \cite{cong2002hybrid,134479}, time of flight (TOF) \cite{cho2010mobile,larsson1996mobile}, received signal strength(RSS) \cite{li2006rss,li2002detection}.

When the source is stationary and the measurements contain no noisy signal, the task can be easily succeeded by getting a small number of measurements. However, in the real world,  these conditions can not be met, therefore, the agent searching for a source requires an estimator to solve the uncertainty issues arisen from the target's motion or the noisy signal which can be studied under adaptive target localization. In \cite{fidan,fidan2013adaptive,fidan2015least}, the authors present a source position estimation algorithm where the agent is able to measure its distance to the position of the source. The algorithm is shown to be exponentially stable under a persistent excitation (PE) condition and robust to drifts in the source location, and the presented simulation results demonstrate that the proposed algorithm performs well in presence of sensor noise as well.

In \cite{fidan2015adaptive}, a geometric cooperative technique is proposed to  estimate permittivity and path loss coefficients for the electromagnetic signal case, with RSS and TOF based range sensors. The proposed technique is integrated to a recursive least squares (RLS)-based adaptive localization scheme and an adaptive motion control law, to perform adaptive target localization robust to uncertainties in environmental signal propagation coefficients. In \cite{umay2017adaptive}, this technique is applied to the problem of tracking biomedical capsule for gastro-intestinal endoscopy and medication applications. 

The above studies all utilize sensor units providing geometric measurements, such as distance, bearing, distance difference, directly related to relative position of the target or the signal source. In many applications, as opposed to distance/direction measurement, RSS is used to estimate the gradient of the unknown signal field of interest and locate the extremum point where the gradient of the field is zero. In \cite{skobeleva2018planar1}, the authors studies a combined formation acquisition and cooperative extremum seeking control scheme for a team of three robots moving a plane in order to find the extremum point of an unknown signal strength field by on-board signal measurement. The proposed algorithm guarantees convergence to a specified neighbourhood of the maximum
of the field while ensuring that the desired formation is acquired and maintained. Similar to the above work, it is accomplished to locate a source by using only direct measurements of that signal at the vehicles’ individual locations in \cite{brinon2016distributed, moore2010source, ogren2004cooperative}.

In this paper, we study adaptive Hessian estimation and extremum localization of a (signal) field $F$ by a sensory agent that continuously measures the intensity of $F$ at its current location while moving. Beyond from the existing literature, including \cite{skobeleva2018planar1, brinon2016distributed, moore2010source, ogren2004cooperative} , the aimed contribution is two-folds: (1) On-line identification of more detailed information about the signal field $F$ than just the extremum of it. (2) More accurate and faster localization of the extremum utilizing this extra information. Having the knowledge of the position $y$ of the sensory agent and the signal value $F(y)$ at the agent's current location as measured by an on-board sensor, we design an adaptive scheme, involving some regression filters, for adaptive estimation of Hessian parameters of $F$, which helps us extract the information of the source location.

Rest of the paper is arranged as follows: The signal map representation is formally introduced and the extremum localization problem is defined in Section \ref{section:The Extremum Localization Problem}. The proposed adaptive Hessian estimation and extremum localization scheme is presented in Section \ref{section:he Proposed Adaptive Hessian Estimation and Localization Scheme}. Stability and the convergence of the proposed scheme are analyzed in Section \ref{section:Stability and Convergence}. Simulation results are displayed to verify the feasibility and robustness of the proposed adaptive scheme in Section \ref{simres}. Concluding remarks are given in Section \ref{sec:Conclusion}.

\section{The Extremum Localization Problem}
\label{section:The Extremum Localization Problem}
The main objective of the adaptive  estimator designs in this paper is to produce an accurate estimate of the location of the extremum(maximum) of a quadratic (signal field) function $ F(\cdot):\mathbb{D} \longrightarrow \mathbb{R}$, for a compact state location domain $\mathbb{D} \subset \mathbb{R}^m$, formulated by

\begin{align}
&&F(y)=c_{1}-\frac{1}{2} \left(y-x\right)^{T}H\left(y-x\right)
\label{eq:1}
\end{align}
where $c_{1}$ is an unknown positive constant and $H$ is an unknown $m {\times} m$ positive definite matrix.  For $m \in \{\ \mathbb{R}^2,\mathbb{R}^3 \}$, \eqref{eq:1} typically represents the strength of a signal emitted by a source at location(state) $x \in \mathbb{R}^m $ measured by a sensory node at location (state) $y \in \mathbb{R}^m$ \cite{krstice1,krstice2,krstice3}. The idea for using a quadratic function as a profile of the signal field is rooted in the fact that any smooth function can be approximated locally by its Taylor expansion near each extremum point. For a general nonlinear smooth function $F_g(\cdot)$, the gradient $\nabla F_g(y)$ will vanish at the extremum point $y = x$, we can write \cite{lang2012calculus} :
\begin{equation}
\label{eq:2}
F_g({x+y_r})=F_g\left(x\right)+\frac{1}{2}y_r^{T}\nabla^{2}F_g(x)y_r+h.o.t
\end{equation}
where $ y_r=y-x$. The approximation (\ref{eq:2}) enables us to extract the gradient of the field using averaging methods \cite{khalil} and find the location of the extremum point. Assuming that $F_g(\cdot)$ is a positive concave signal field function, $\nabla^2F_g(x)$ is negative definite and $c_1$ and $H$ in \eqref{eq:1} matches, respectively, with $F_g(x)$ and $-\nabla^2F_g(x)$ in \eqref{eq:2}. For brevity, neglecting the higher order terms ( $h.o.t.$ ) in \eqref{eq:2}, we focus on the representation \eqref{eq:1} in this paper, and formally define the extremum localization problem for this representation.

\par\textit{Problem 1}: Consider the quadratic signal field function in \eqref{eq:1}. Suppose that a sensory agent has access to the field measurement $F(y)$ at its current location $y$. Design an adaptive identification scheme to estimate the target location $x$ at which $F$ takes its maximum value, and derive the conditions under which the estimate $\hat x(t)$ converges to $x$ asymptotically.

\section{The Proposed Adaptive Hessian Estimation and Localization Scheme}
\label{section:he Proposed Adaptive Hessian Estimation and Localization Scheme}
In order to devise an adaptive localization algorithm, we use the adaptive parameter identification based framework proposed in \cite{fidan,fidan2013adaptive,fidan2015least}. We use the notation in \cite{fidan} for derivative operation and asymptotically equal signals: $s$ denotes the derivative operator, i.e., given a function $f$ of time $t$,  
$s f := \dot{f}=d f /dt$. $\frac{1}{s+a}f(t):=\int_0^t{e^{-a\tau}f{\tau}d\tau}$. For two vector functions $f,g$ of the same dimension, $f(\cdot)\approx g(\cdot)$ if there exist $\lambda, M$ such that 
$\|f(t)-g(t)\| \leq M e^{-\lambda t}$ for all $t\geq 0$. We derive a parametric model that is linear in unknown parameters of the system, i.e., the elements of Hessian matrix $H$ and the location(state) $x$ of the extremum. Taking time derivative of \eqref{eq:1} and assuming that $x$ is constant,
i.e., $\dot{x}=0$,
we obtain
\begin{align}
\dot F (y) =&- \dot y^{T}H (y-x) = - \dot y^{T}H y+\dot y^{T}H x \nonumber \\
=&-\frac{1}{2} \frac{d}{dt}\left(y^THy \right)+\frac{d}{dt}\left(y^T\right)Hx \nonumber \\
=& -\frac{1}{2}\frac{d}{dt} \big ( H_{11}y_1^2+2H_{12}y_1y_2+\cdots+H_{22}y_2^2 \nonumber \\
&+2H_{23}y_2y_3+\cdots+H_{mm}y_m^2 \big)+\frac{d}{dt}\left(y^T\right)Hx
\label{eq:4}
\end{align}
which can be written as 
\begin{align}
&& \dot F (y)= \theta^{*T} \frac{d\Psi}{dt},
\label{2DidotDi}
\end{align}

\begin{equation}
\begin{array}{ll}
\label{eq:5}
\theta^{\ast}=\bigg[&H_{11},H_{12},
 \cdots,H_{1m},H_{22},\cdots,H_{mm}, \\ & 
\underbrace{x^TH_1,\cdots,x^TH_m}_{x^TH} \bigg]^{T}\in\mathbb{R}^{\frac{m(m+3)}{2}},
\end{array}
\end{equation}

\begin{equation}
\begin{array}{ll}
\Psi = \bigg  [ &\frac{-1}{2} y_{1}^2 ,-y_{1}y_2 , \cdots , -y_{1}y_{m} ,\frac{-1}{2} y_{2}^2 , \\ &  \cdots ,  \frac{-1}{2}y_{m}^2 ,  y^T \bigg ]^{T}  \in\mathbb{R}^{\frac{m(m+3)}{2}},
\label{eq:6m}
\end{array}
\end{equation}
where $H_i$ denotes the $i$th column (= transpose of the $i$th row) of $H$. In order to eliminate need for explicit differentiation of available signals, $z(\cdot)$ and $\phi(\cdot)$ are introduced as the state variable filtered versions of $F(\cdot)$ and $\Psi(\cdot)$, respectively:
\begin{align}
\label{eqn:filterz:single:start}
   && \dot \xi_1(t)&=-a \xi_1(t) + F(y(t)), \\ && \xi_1(0)&=0,   \\ &&
    \label{eqn:filterz:single:mid}
    z(t)&=-a\xi_1(t)+F(y(t)), \\ &&
    \dot \xi_2(t)&=-a \xi_2(t) + \Psi(t), \\ &&\xi_2(0)&=[0,\dots,0]^T\in \mathbb{R}^{\frac{m(m+3)}{2}}, \\ &&
    \label{eqn:filterz:single:end}
    \phi(t)&=-a\xi_2(t)+\Psi(t), 
\end{align}
for some $a>0$. It can be seen in  \eqref{eqn:filterz:single:start}--\eqref{eqn:filterz:single:end} that the measurements of the location(state) $y(t)$ of the sensory agent and the field intensity $F(y(t))$ at that location are sufficient to generate the signals $z(t)$ and $\phi(t)$. 

\begin{lem}
\label{lemma:1}
Suppose $\theta^* \in \mathbb{R}^{\frac{m(m+3)}{2}}$ is a constant, and $z(t), \phi(t)$ are defined by \eqref{eqn:filterz:single:start}--\eqref{eqn:filterz:single:end} with $a>0$. Then there holds:
\begin{align}
\label{lemmaeq1}
&&z(\cdot) \approx \theta^{*T} \phi(\cdot).
\end{align}
\end{lem}

\par\textit{Proof}:Using \eqref{eqn:filterz:single:start}--\eqref{eqn:filterz:single:mid},  we  obtain;
\begin{align}
&&\dot z (t) + a z(t)=\frac{d}{dt} \left\{ F \right \},
\end{align}
where $a>0$. In operator notation i.e., using $s$ to denote the differentiator operator,
\begin{align}
&&z(\cdot) \approx \frac{s}{s+a}  \left\{ F(\cdot)  \right \}.
\end{align}
Similarly,
\begin{align}
&&\phi(\cdot)  \approx \frac{s}{s+a}  \left\{ \Psi(\cdot)  \right \}.
\end{align}
Then, 
\begin{align}
z(\cdot) \approx \frac{s}{s+a}   \left\{ F(\cdot)  \right \} \approx  \frac{1}{s+a}   \left\{ \theta^{*T} \dot \Psi(\cdot)  \right \} \nonumber \\ \approx 
\theta^{*T} \frac{s}{s+a}   \left\{ \Psi(\cdot)  \right \} \approx \theta^{*T} \phi(\cdot). \  \blacksquare
\end{align}

Using  \eqref{lemmaeq1} as linear parametric model, and \eqref{eqn:filterz:single:start}--\eqref{eqn:filterz:single:end}  to generate the regressor signals in this model, we design the following gradient based adaptive estimation algorithm \cite{ion,ioannou2006adaptive} to identify $\theta^*$:

\begin{align}
%\dot{\hat{\eta}}=\gamma\phi\left(z-\hat{\eta}\right)=\gamma\phi\phi^{T}(\eta^{\ast}-\hat{\eta}), \ \ \tilde{\eta}=-\eta^{\ast}+\hat{\eta}
&&\dot {\hat \theta} = \gamma \phi (z- \hat \theta^T \phi),
\label{eq:8}
\end{align} 
 where $\hat{\theta}$ denotes the estimate of ${\theta}^*$ and $ \gamma > 0$ is a scalar design constant. To be able to extract the information of the elements of $H$ and the location(state) of the source ($x$) from the estimation of $ \theta^*$, we consider the following partitioning of $ \theta^*$ and  $\hat \theta$ ;
 
\begin{align}
&& \theta^*= \begin{bmatrix}
 \theta_H^* \\
\theta_x^*
\end{bmatrix}, \ \
\hat \theta= \begin{bmatrix}
\hat \theta_H \\
\hat \theta_x
\end{bmatrix}
\end{align}
where $\theta_H^* \in \mathbb{R}^{\frac{m(m+1)}{2}} $ is composed of the entries of $\theta^*$ that are independent of $x$, $\theta_x^*= H  x \in \mathbb{R}^m $, $\hat \theta_H$ and $\hat \theta$ are the estimates of $\theta_H^*$ and $\theta_x^*$ respectively. Since all the elements of $H$ exist in $\theta_H^*$, we can form $\hat H$ (the estimate of $H$) from $\hat \theta_H$. In order to obtain $\hat x$ which is the estimation of the source's location(state) $x$, we utilize the equality $\theta_x^*=Hx$;
\begin{align}
\label{estimationofsource}
&&\hat x= \hat H^{-1} \hat \theta_x.
\end{align}
In order to take the inverse of $\hat H $ in \eqref{estimationofsource}, it must be guaranteed that $\hat H$ is non-singular. 

\begin{assumption}
\label{assump:Hpositive}
The Hermitian matrix $H$ satisfies the following: 
\begin{enumerate}
    \item $H_{ii}>0$ for all $i=1,\cdots,m$.
    \item $H$ is strictly diagonally dominant which means $|H_{ii}|>\sum_{i \neq j}|H_{ij}|$ for all $i,j=1,\cdots,m$.
\end{enumerate}
\end{assumption}
\begin{lem}
If $H$ satisfies Assumption \ref{assump:Hpositive}, then it is positive definite.
\end{lem}
\par\textit{Proof}: The result is a direct corollary of Theorem 6.1.10 of \cite{horn1990matrix}.

%Since $v_i^THv_i>0$ for any $v_i \in \mathbb{R}^m$, we can prove $H_{ii}>0$ by taking $v_i$ as a vector of all zeros, except for a 1 in the $i^{th}$ place. \blacksquare

To assure $\hat H$ is non-singular, we apply parameter projection on the elements of $\hat \theta_H$ in consideration of Assumption \ref{assump:Hpositive} and \eqref{eq:8} with the parameter projection is re-designed as;
\begin{align}
\label{eq:8proj}
&&\dot{\hat \theta}= \underset{\hat \theta_H \in S_H }{\text{Proj} } \{ \gamma \phi (z- \hat \theta^T \phi) \},
\end{align}
where the convex compact set $S_H$ is defined as the set of all vectors $\hat \theta_H=[\hat H_{11},\hat H_{12},\cdots,\hat H_{1m},\hat H_{22},\cdots,\hat H_{mm} ]^T$ such that the corresponding $m \times m$ matrix $\hat H$ satisfies Assumption \ref{assump:Hpositive}, and $\underset{\hat \theta_H \in S_H }{\text{Proj} } \{\cdot \}$ is the parameter projection operator \cite{ion,ioannou2006adaptive} defined to maintain $\hat \theta_H$ in $S_H$.

\emph{Remark 2.1} If $H$ is a diagonal matrix, the vectors $\theta^*$ and $\Psi$ in \eqref{eq:5}--\eqref{eq:6m} can be redefined in reduced form as follows:
\begin{align}
\label{reducedtheta}
&&\theta^{\ast}=&\bigg[H_{11},\cdots, H_{mm}, x^TH\bigg]^{T}\in\mathbb{R}^{2m} \\ 
&&  
\label{reducedphi} \Psi=&\bigg[\frac{-1}{2}y_{1}^2,\cdots,\frac{-1}{2}y_{m}^2,
 y^T\bigg]^{T}\in\mathbb{R}^{2m}  
\end{align}

For a general case, since $H$ is a symmetric matrix with real elements, we can deduce that by choosing appropriate coordinates, we can diagonalize the matrix $H$ and hence, design the identification algorithm  based on the reduced order model \eqref{lemmaeq1},\eqref{reducedtheta},\eqref{reducedphi}.

In the next section, we analyze the stability of the proposed adaptive estimation and localization scheme.

\section{Stability and Convergence}
\label{section:Stability and Convergence}

\subsection{Stationary Extremum Localization}

Note that the base adaptive law \eqref{eq:8} and the adaptive law \eqref{eq:8proj} with parameter projections can be rewritten, respectively, as
\begin{align}
\label{etaconvergesnew}
&&\dot{\tilde \theta}=& \dot{\hat \theta} =  - \gamma \phi \phi^T \tilde \theta , \\
\label{etaconverges}
&&\dot{\tilde \theta}=& \dot{\hat \theta} =  \underset{\hat \theta_H \in S_H }{\text{Proj} } \{- \gamma \phi \phi^T \tilde \theta \}, 
\end{align} 
where $\tilde \theta = \hat \theta - \theta^*$.
Hence, the aimed convergence of the estimate $\hat \theta$ to actual $\theta^*$ is equivalent to the convergence of $\tilde \theta$ to zero. 

\begin{thm}
\label{thm:localization}
Suppose $\theta^* \in \mathbb{R}^{\frac{m(m+3)}{2}}$ is a constant. Consider $z(t)$ and $\phi(t)$ defined in \eqref{eqn:filterz:single:start}--\eqref{eqn:filterz:single:end}, with $a>0$. Then for each of the base adaptive law \eqref{etaconvergesnew} and the adaptive law \eqref{etaconverges} with parameter projection , there exist $\rho_1, \rho_2, \lambda>0$ such that for all $t \geq 0$ and $||\theta^*(0)||$
\begin{align}
&&||\tilde \theta (t) || \leq (\rho_1||\theta^*(0)|| + \rho_2 )e^{-\lambda t} 
\end{align}
if and only if there exist $\alpha_1 >0$, $\alpha_2>0$, $T>0$ such that for all $t \geq 0$
\begin{align}
\label{boundsoftheorem2}
&&\alpha_1I \leq \int_{t}^{t+T}  \phi(\tau) \phi(\tau)^T d\tau \leq \alpha_2 I.
\end{align}
\end{thm}
\par\textit{Proof}: It is established in the literature (see, e.g., \cite{anderson1977exponential}) that \eqref{etaconvergesnew} is exponentially asymptotically stable if and only if \eqref{boundsoftheorem2} holds. Moreover, it is proven in \cite{ioannou2006adaptive} that the parameter projection does not affect the properties of the gradient adaptive laws deducted on the Lyapunov analysis and it can only make the time derivative of Lyapunov function more negative. Hence, \eqref{etaconverges} is also exponentially asymptotically stable if and only if \eqref{boundsoftheorem2} holds.$\blacksquare$

\subsection{Drift in Extremum Location}
\label{subsec:driftin}
The drift analysis in \cite{fidan} can be applied here as well, without requiring significant modification.
Before, detailing the drift analysis, we make the following assumption.
\begin{assumption}
\label{boundsofyandx}
The agent trajectory $y:\mathbb{R} \rightarrow \mathbb{R}^m $ is twice differentiable, the source trajectory $x:\mathbb{R} \rightarrow \mathbb{R}^m $ is differentiable and there exist $M_1,M_2,M_3,M_4,\epsilon>0$ such that for all $t \in \mathbb{R}$
\begin{align}
&&||y(t)||& \leq M_1, \ || \dot y(t)|| \leq M_2, \ || \ddot y(t)|| \leq M_3,\\
   &&||x(t)|| &\leq M_4, \ \label{xdotbound}
    ||\dot x(t)|| \leq \epsilon.
\end{align}
\end{assumption}

\begin{lem}
\label{lem:driftstart}
Under Assumption \ref{boundsofyandx}, for $z(t) $ and $\phi(t)$ defined in \eqref{eqn:filterz:single:start}--\eqref{eqn:filterz:single:end}, there exists $M_5 :\mathbb{R}_{\geq 0} \rightarrow \mathbb{R}_{\geq 0}$ such that for a suitable $K_1$ depending only on $M_1,M_2,M_4$ and $a$,
\begin{align}
\label{lem:bounds}
  &&  |z(t)-\theta^{*T}\phi(t)| \leq M_5(t), \ \ \ \forall t \geq 0
\end{align}
and
\begin{align}
&&M_5(\cdot) \approx K_1 \epsilon   .
\end{align}

\end{lem}
\par\textit{Proof}:
 Using the operator notation in the proof of Lemma \ref{lemma:1}, it is achieved that
\begin{align}
\label{zdrift}
    z(\cdot) \approx &\frac{s}{s+a}\left\{F(\cdot)\right\} \nonumber \\ 
    \approx &\frac{1}{s+a} \left \{-\left(\dot y(\cdot) - \dot x(\cdot) \right)^T H \left( y(\cdot) -  x(\cdot) \right) \right \} \nonumber \\ 
     \approx & - \frac{s}{s+a}\left \{ \frac{1}{2} y^T(\cdot)Hy(\cdot)\right \} \nonumber \\ 
    &+ \frac{1}{s+a}\left \{ \frac{1}{2} \dot y^T(\cdot)Hx(\cdot)\right \} + f(\cdot)
\end{align}
where 
\begin{align}
   && f(\cdot)= \frac{1}{s+a}\left \{ \frac{1}{2} \dot x^T(\cdot)H \left(y(\cdot)- x(\cdot)\right) \right \}.
\end{align}

In consideration of Assumption \ref{boundsofyandx}, there exists a $F;\mathbb{R}_{\geq 0} \rightarrow \mathbb{R}_{\geq 0}$, such that for all $t \geq 0$,
\begin{align}
  &&  |f(t)| \leq F(t)
\end{align}
and
\begin{align}
  &&  F(\cdot) \approx \frac{M_1+M_4}{a}\epsilon.
\end{align}

Now, consider the second term in \eqref{zdrift}
\begin{align}
 &&   \frac{1}{s+a}\left \{ \frac{1}{2} \dot y^T(\cdot)Hx(\cdot)\right \} \approx Q(\cdot)
\end{align}
where with $C \in \mathbb{R}^m$,
\begin{align}
    Q(t)=&e^{-at} \int_0^t e^{a\tau} \dot y^T(\tau) H x(\tau) d\tau \nonumber \\
    =&e^{-at} \left [ \left ( \int_0^\tau e^{as} \dot y (s) ds +C  \right)^T Hx(\tau) \right]_0^t \nonumber \\
    &-e^{-at} \int_0^t \left( \int_0^\tau e^{as} \dot y (s) ds +C\right)^T H \dot x(\tau)d \tau \nonumber \\
    =& \left [ \left ( \int_0^\tau e^{-a(t-s)} \dot y (s) ds +C e^{-at} \right)^T Hx(\tau) \right]_0^t\nonumber \\
    &-G(t), \\
\label{Gt}
    G(t)=&e^{-at} \int_0^t \left( \int_0^\tau e^{as} \dot y (s) ds +C\right)^T H \dot x(\tau)d \tau .
\end{align}
Thus, as $a>0$, and adding the first term in \eqref{zdrift}, we obtain
\begin{align}
- \frac{s}{s+a}\left \{ \frac{1}{2} y^T(\cdot)Hy(\cdot)\right \}+Q(\cdot) \approx {\theta^*}^T \phi (\cdot)-G(\cdot).
\end{align}
Moreover, from \eqref{Gt}, it is obtained that 
\begin{align}
\label{Gtbounds}
    |G(t)| \leq e^{-at}M_2 \lambda_{max}(H) \epsilon \left[ \frac{e^{at}-1}{a^2}+t \left( ||C||-\frac{1}{a} \right) \right ]
\end{align}
Then the result follows from \eqref{lem:bounds}--\eqref{Gtbounds}. $\blacksquare$

%Then, \eqref{zdrift} becomes
%\begin{align}
%    z(\cdot) \approx Q(\cdot) - G(\cdot) + f(\cdot)
%\end{align}

Then in the view of Theorem \ref{thm:localization}, we have the following result.

\begin{thm}
\label{thm:drift}
Suppose Assumption \ref{boundsofyandx} hold, and there exist $\alpha_1,\alpha_2,T > 0$ such that $\forall t \geq 0$. Consider $z(t)$ and $\phi(t)$ defined in \eqref{eqn:filterz:single:start}--\eqref{eqn:filterz:single:end}. Then $\hat \theta(t)$ in \eqref{eq:8proj} obeys for some K obtained from $M_1,M_2,M_4,\gamma,a,T,\alpha_1$ and $\alpha_2$,  $\text{lim } \text{sup}_{t \rightarrow \infty}|\hat \theta(t)-\theta^*(t)|=K\epsilon$.
\end{thm}
\par\textit{Proof}: Due to \eqref{eq:8proj} there holds
\begin{align}
    \dot{\tilde \theta}(t)=&\dot{\hat \theta}(t)-\dot{ \theta}^*(t)\nonumber \\
    =&\gamma \phi(t) (z(t)- \hat \theta^T(t) \phi(t))-\dot{ \theta}^*(t)\nonumber  \\
    =&-\gamma \phi(t) \phi^T(t)\tilde \theta(t)+\gamma \phi(t) (z(t)-  \theta^{*T}(t) \phi(t))\nonumber  \\
    &-\dot{ \theta}^*(t) \nonumber \\
    =&-\gamma \phi(t) \phi^T(t)\tilde \theta(t)+G_2(t)
\end{align}
where
\begin{align}
   && G_2(t)=\gamma \phi(t) (z(t)-  \theta^{*T}(t) \phi(t))-\dot{ \theta}^*(t).
\end{align}

Then because of Lemma \ref{lem:driftstart}, \eqref{xdotbound} and the fact that $\hat \phi (\cdot)$ is bounded, there exists a $K_5>0$ obtained from $M_1,M_2M_4,\gamma$ and $a$, and an $M_6 : \mathbb{R}_{\geq 0} \rightarrow \mathbb{R}_{\geq 0}$, obeying $M_6(\cdot)\approx K_5 \epsilon$ such that $|G_2(t)| \leq M_6(t) \forall t \geq 0$. Hence the result follows from the exponential asymptotic stability of \eqref{eq:8proj}. $\blacksquare$

\section{Simulation Results}
\label{simres}

In this section, we provide simulation results to exhibit the performance of the proposed scheme in Section \ref{section:he Proposed Adaptive Hessian Estimation and Localization Scheme}. For all examples, the state number, the adaptation gain and the filter pole are selected as $m=2$(considering the localization of extremum in 2-D plane.), $\gamma=1$ and $a=0.5$, respectively and the signal field is formed as $F(y)=3-(y-x)H(y-x)$ where the Hessian matrix is $H=\begin{bmatrix}1&0.2\\0.2&2\end{bmatrix}$.

\begin{figure*}
\begin{subfigure}[t]{.33\linewidth}
\includegraphics[width=\linewidth]{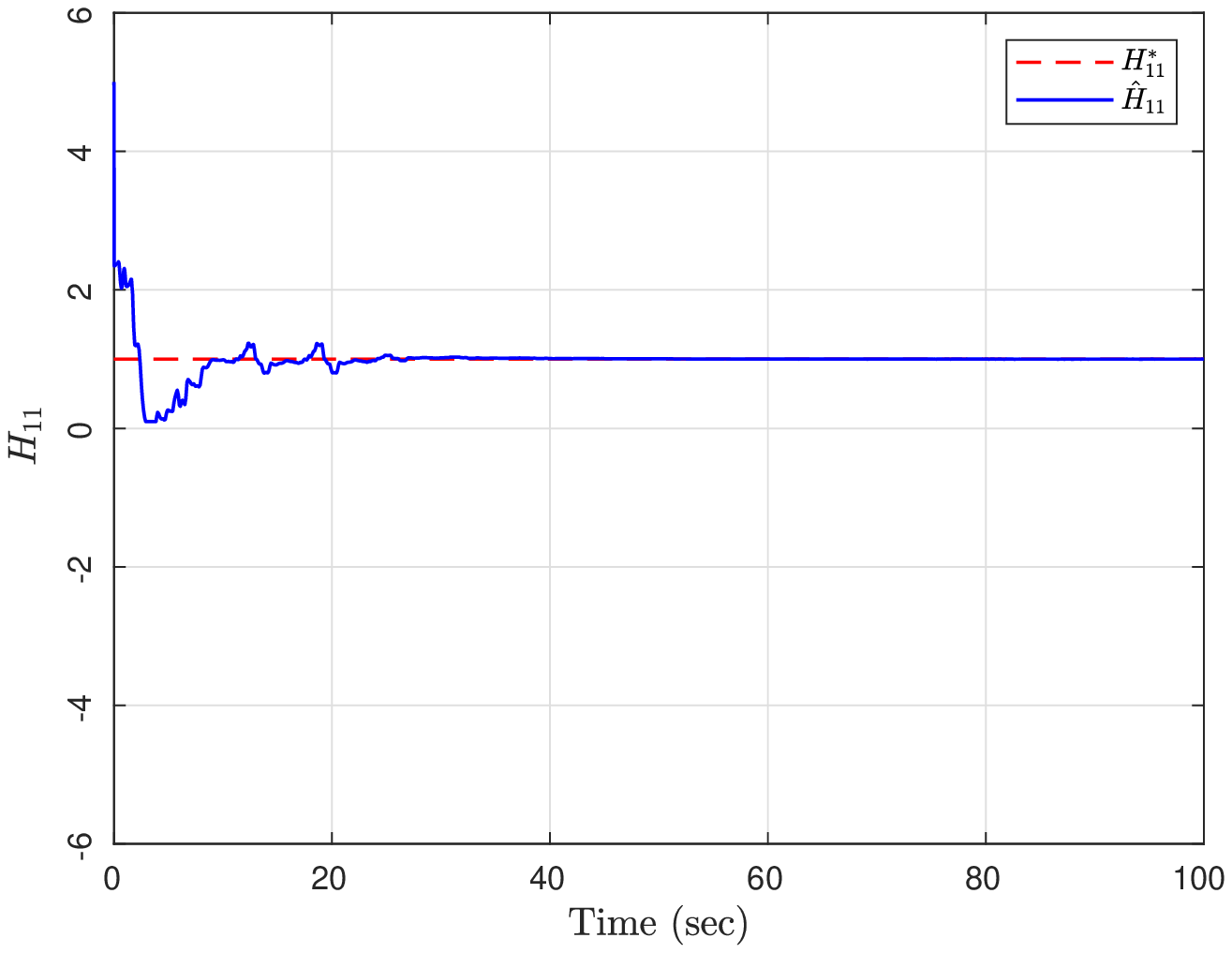}
\end{subfigure}
\begin{subfigure}[t]{.33\linewidth}
\includegraphics[width=\linewidth]{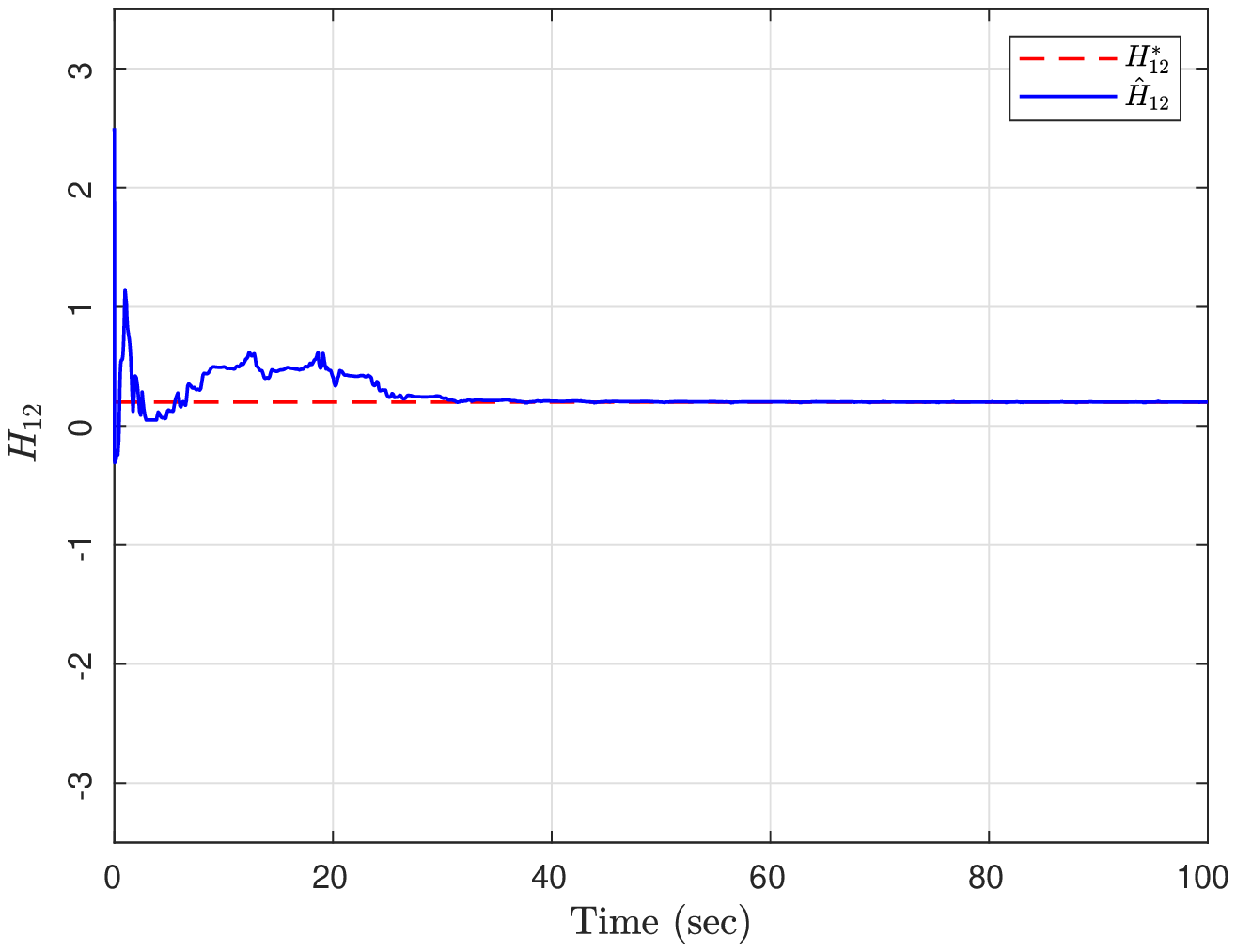}
\end{subfigure}
\begin{subfigure}[t]{.33\linewidth}
\includegraphics[width=\linewidth]{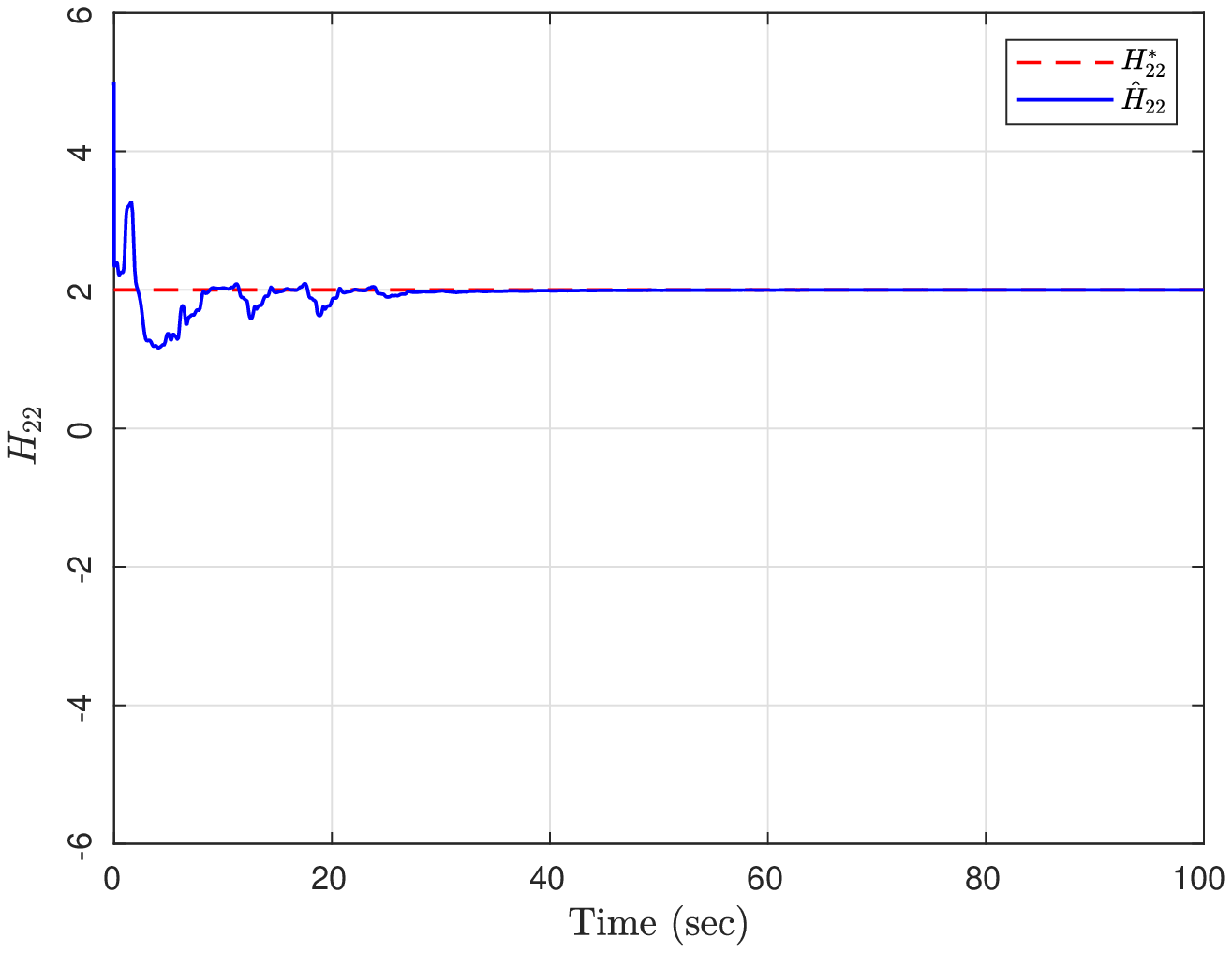}
\end{subfigure}
\begin{subfigure}[t]{.33\linewidth}
\includegraphics[width=\linewidth]{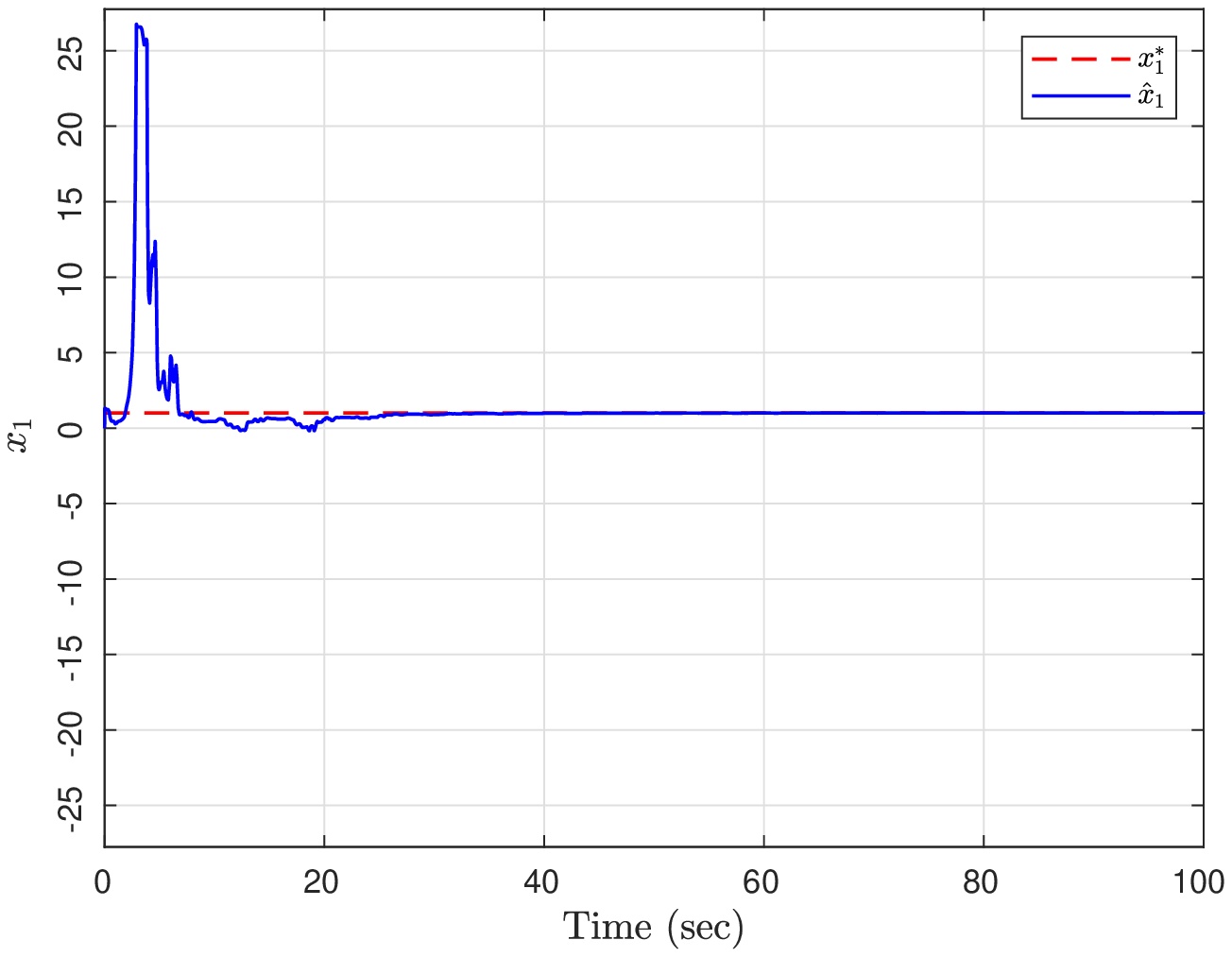}
\end{subfigure}
\begin{subfigure}[t]{.33\linewidth}
\includegraphics[width=\linewidth]{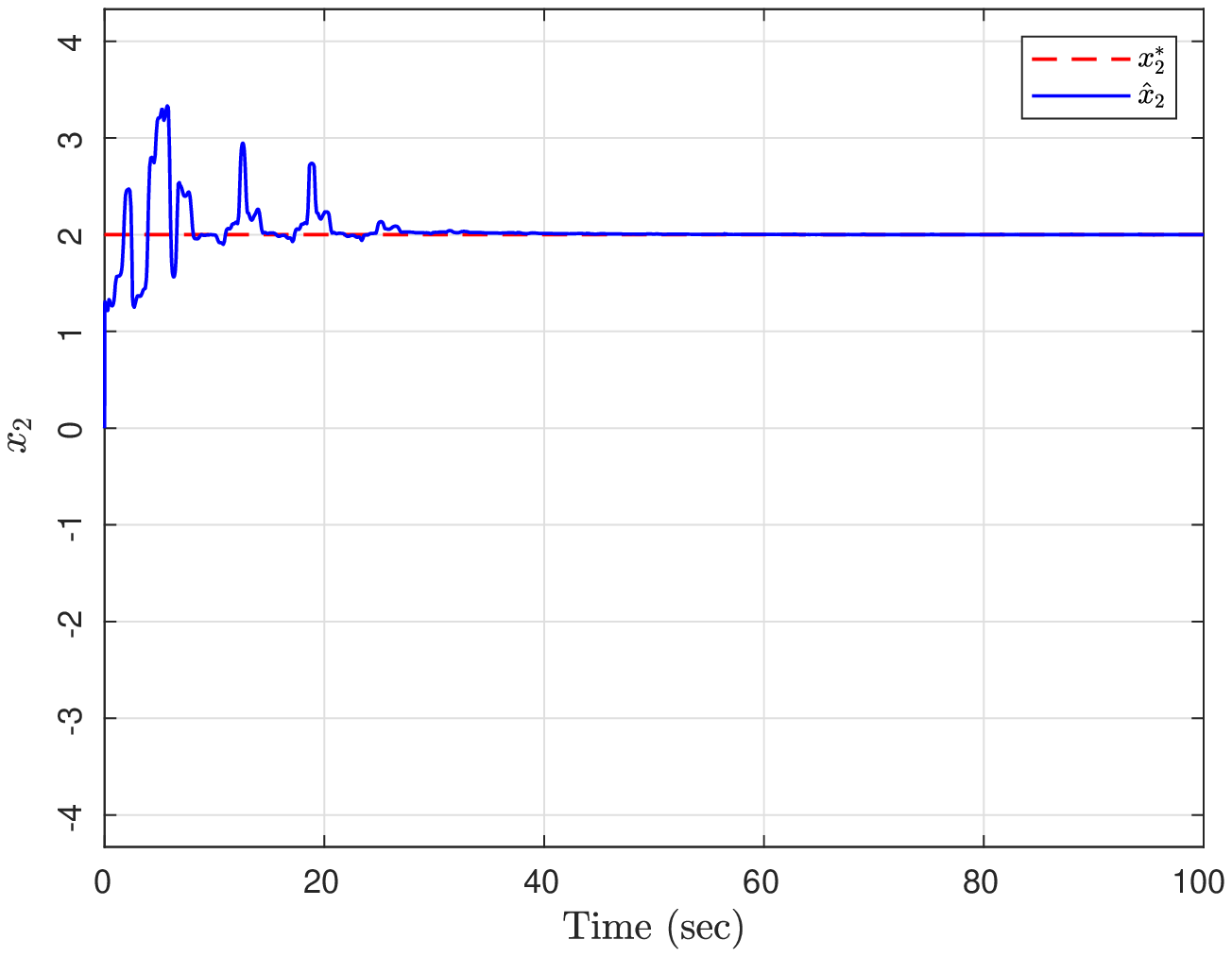}
\end{subfigure}
\begin{subfigure}[t]{.33\linewidth}
\includegraphics[width=\linewidth]{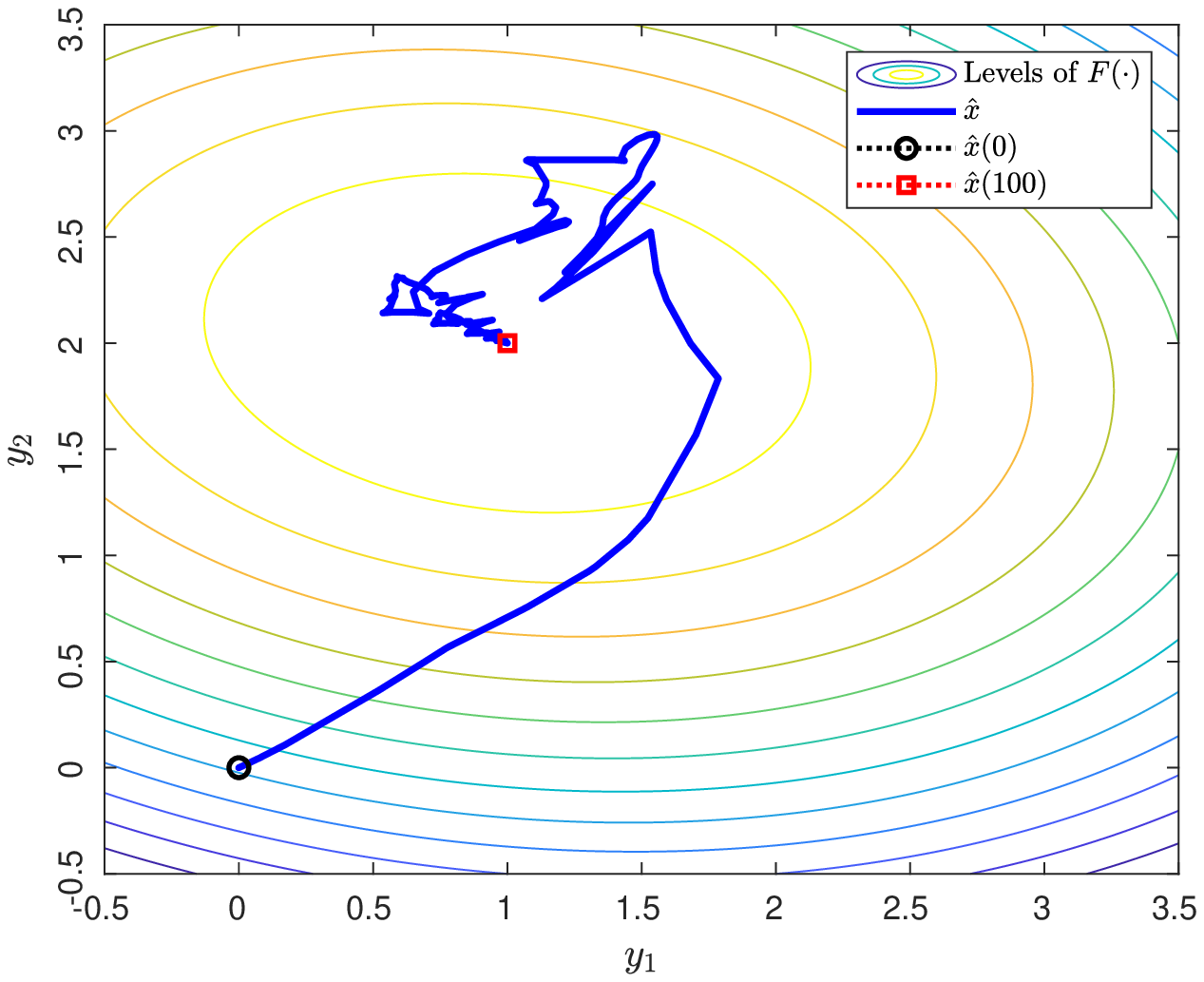}
\end{subfigure}
\caption{ Location estimation for $x(t)=[1, \ 2]^T$, $y(t)=[\sin(4t)+\sin(5t), \ \sin(2t)+\sin(3t)]^T$, $a=0.5$. The dashed lines and the solid lines represent the actual values and  their estimates, respectively.}
\label{fig:nn_nd}
\end{figure*}

Scenario 1: Assume the extremum location is at $x=\begin{bmatrix}1&2\end{bmatrix}^T$ and the sensory agent's trajectory is given by $y = \begin{bmatrix} \sin(4t) + \sin(5t) & \sin(2t) + \sin(3t) \end{bmatrix}^T$. Using the adaptive estimation algorithm \eqref{eq:8proj}, the Hessian matrix and the source location estimates converge to their actual values exponentially as seen in Figure \ref{fig:nn_nd}.

\begin{figure*}
\begin{subfigure}[t]{.33\linewidth}
\includegraphics[width=\linewidth]{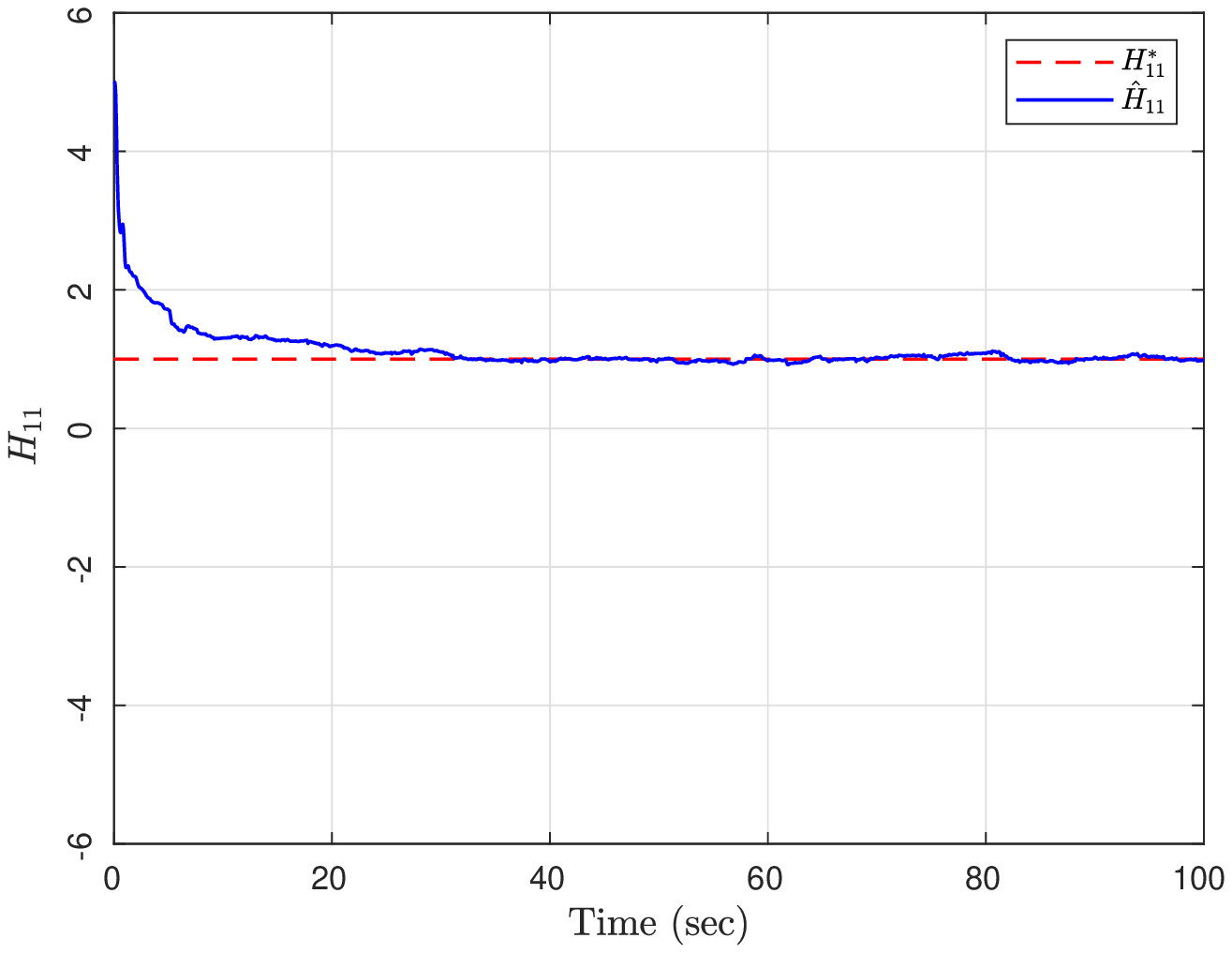}
\end{subfigure}
\begin{subfigure}[t]{.33\linewidth}
\includegraphics[width=\linewidth]{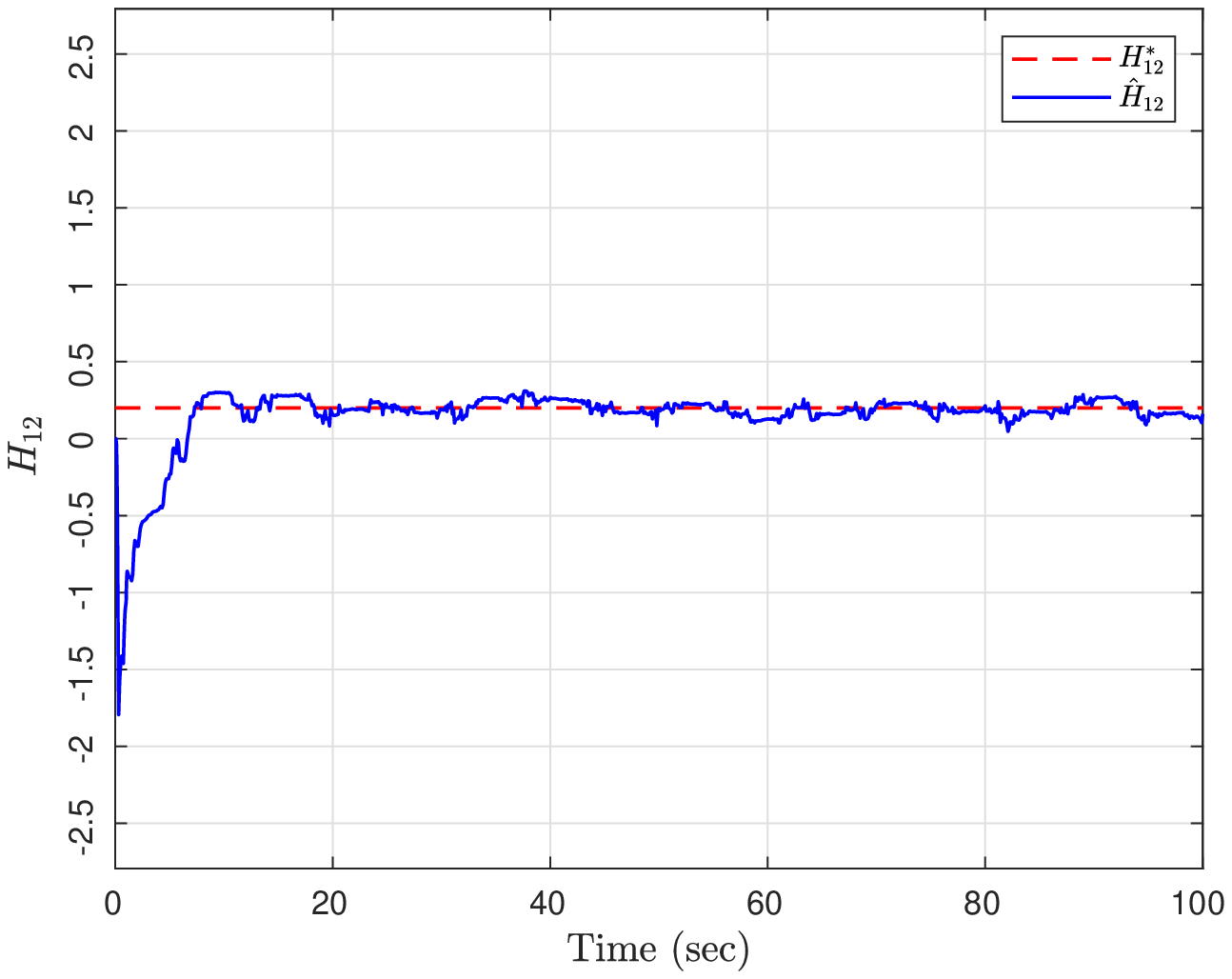}
\end{subfigure}
\begin{subfigure}[t]{.33\linewidth}
\includegraphics[width=\linewidth]{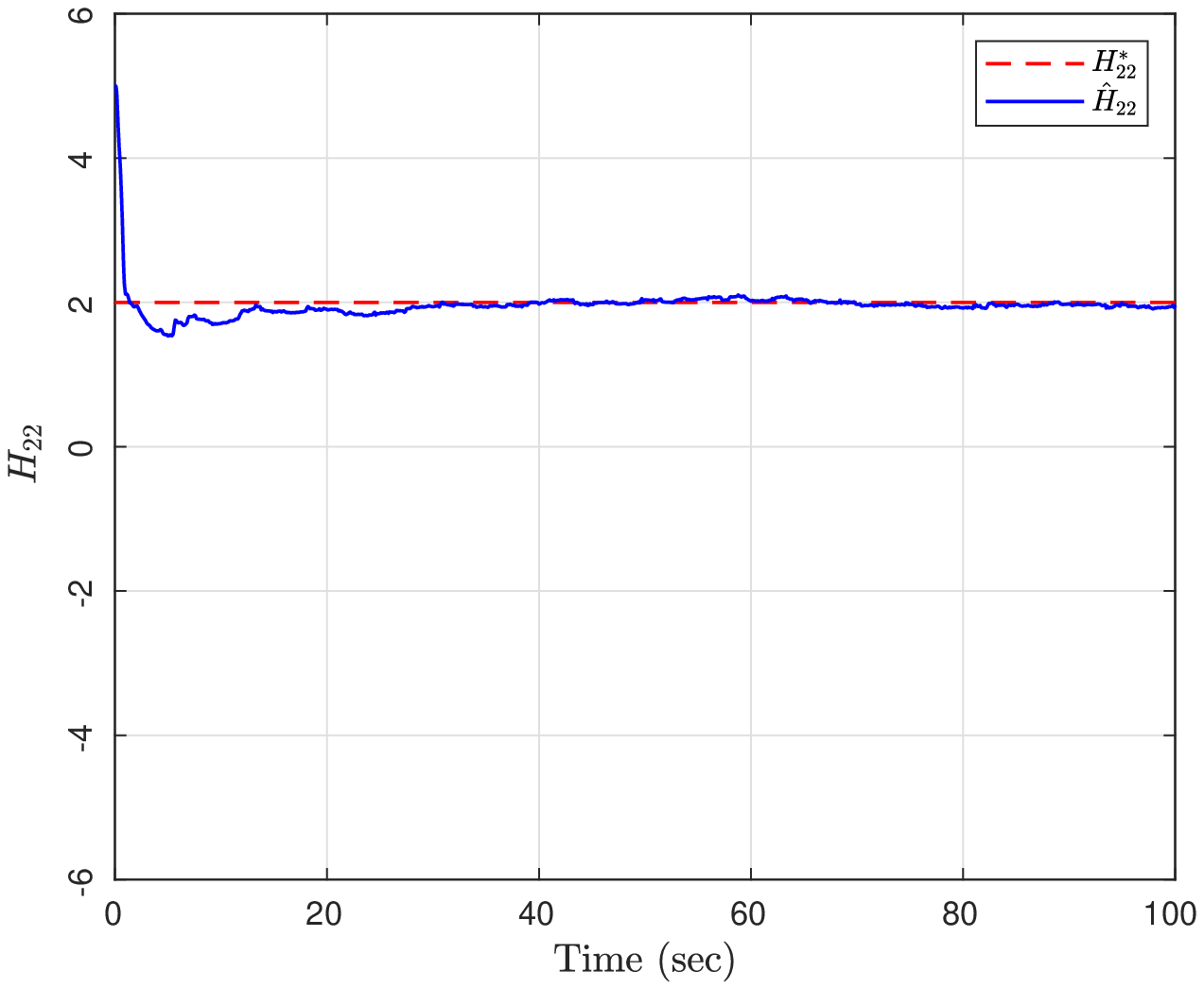}
\end{subfigure}
\begin{subfigure}[t]{.33\linewidth}
\includegraphics[width=\linewidth]{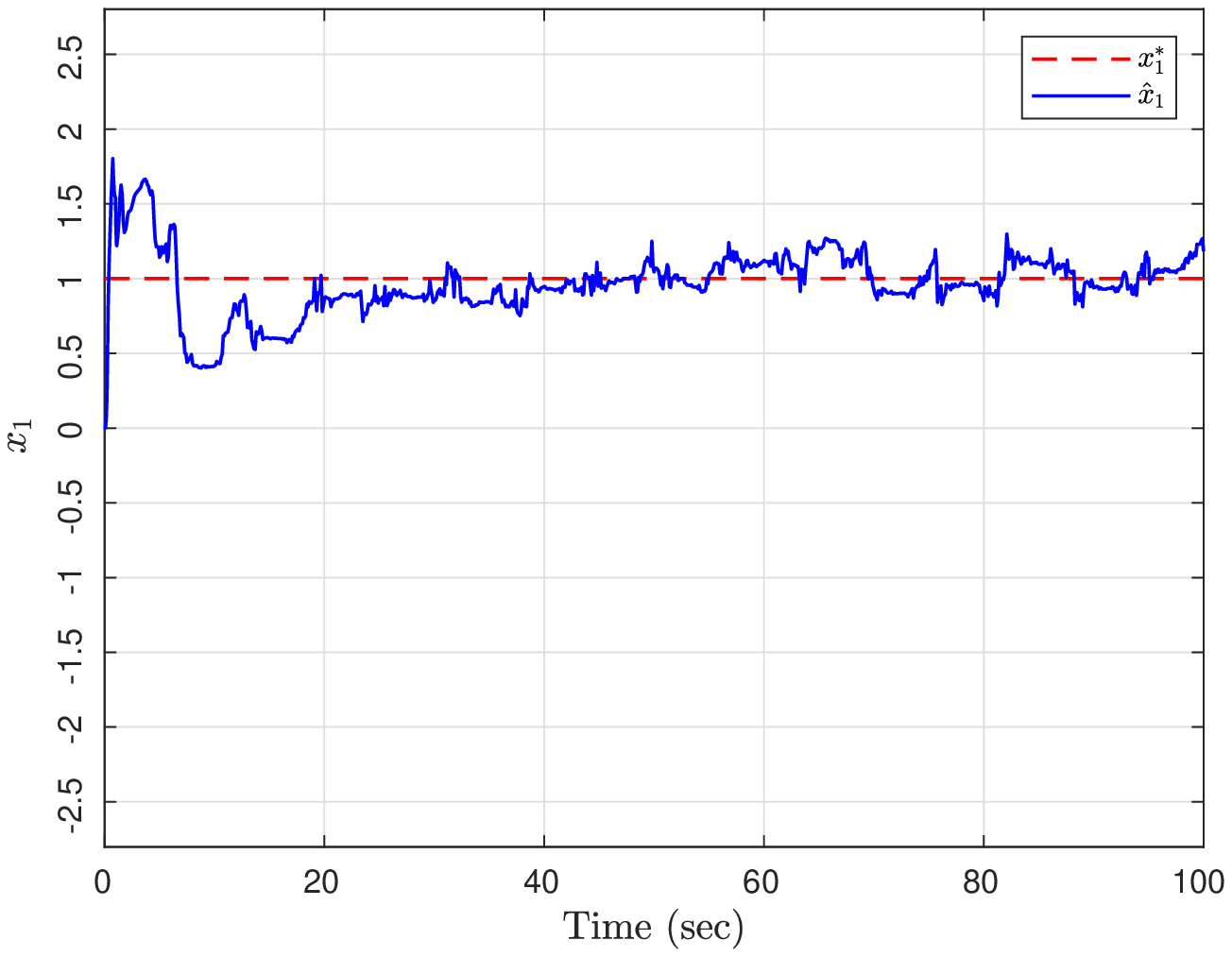}
\end{subfigure}
\begin{subfigure}[t]{.33\linewidth}
\includegraphics[width=\linewidth]{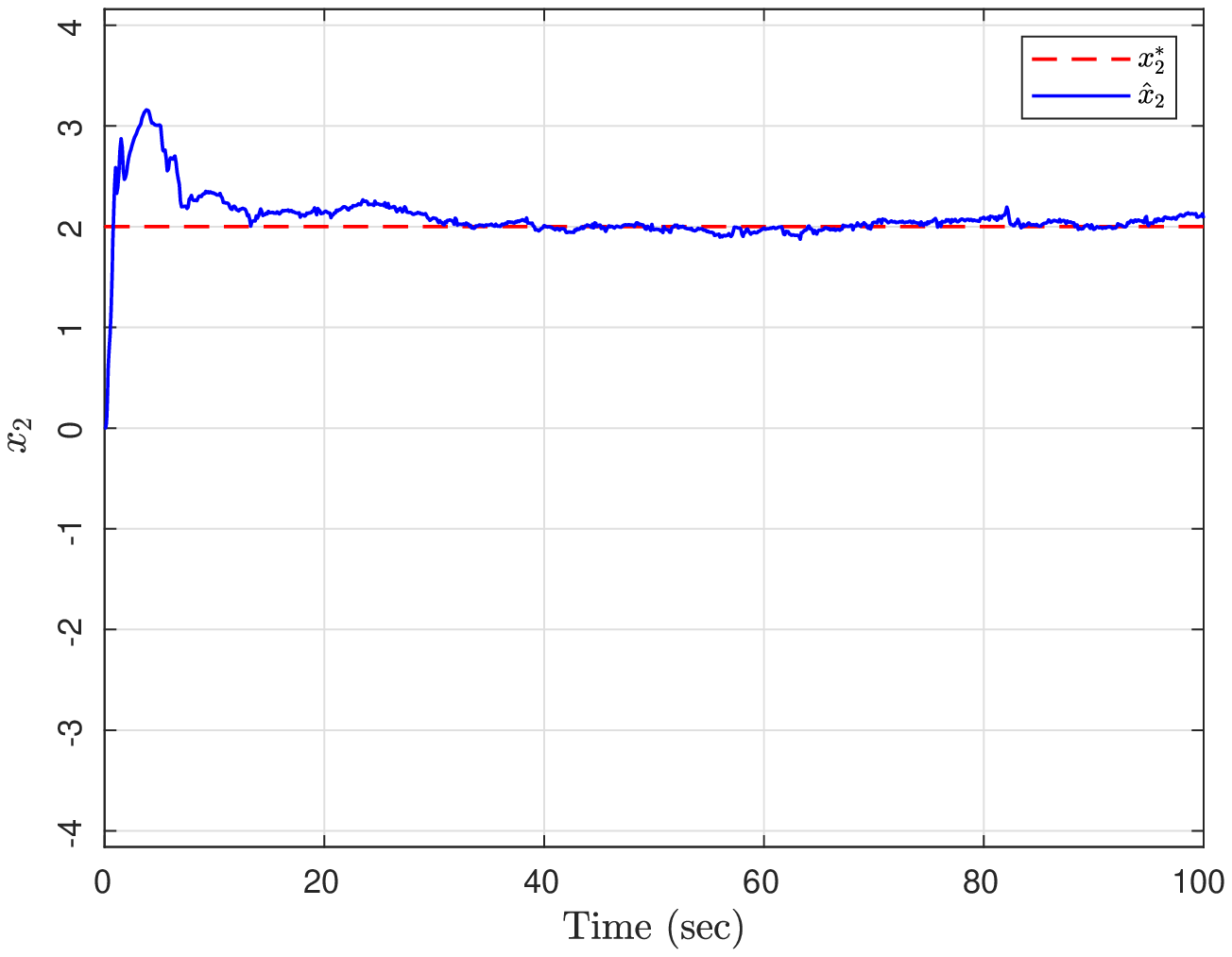}
\end{subfigure}
\begin{subfigure}[t]{.33\linewidth}
\includegraphics[width=\linewidth]{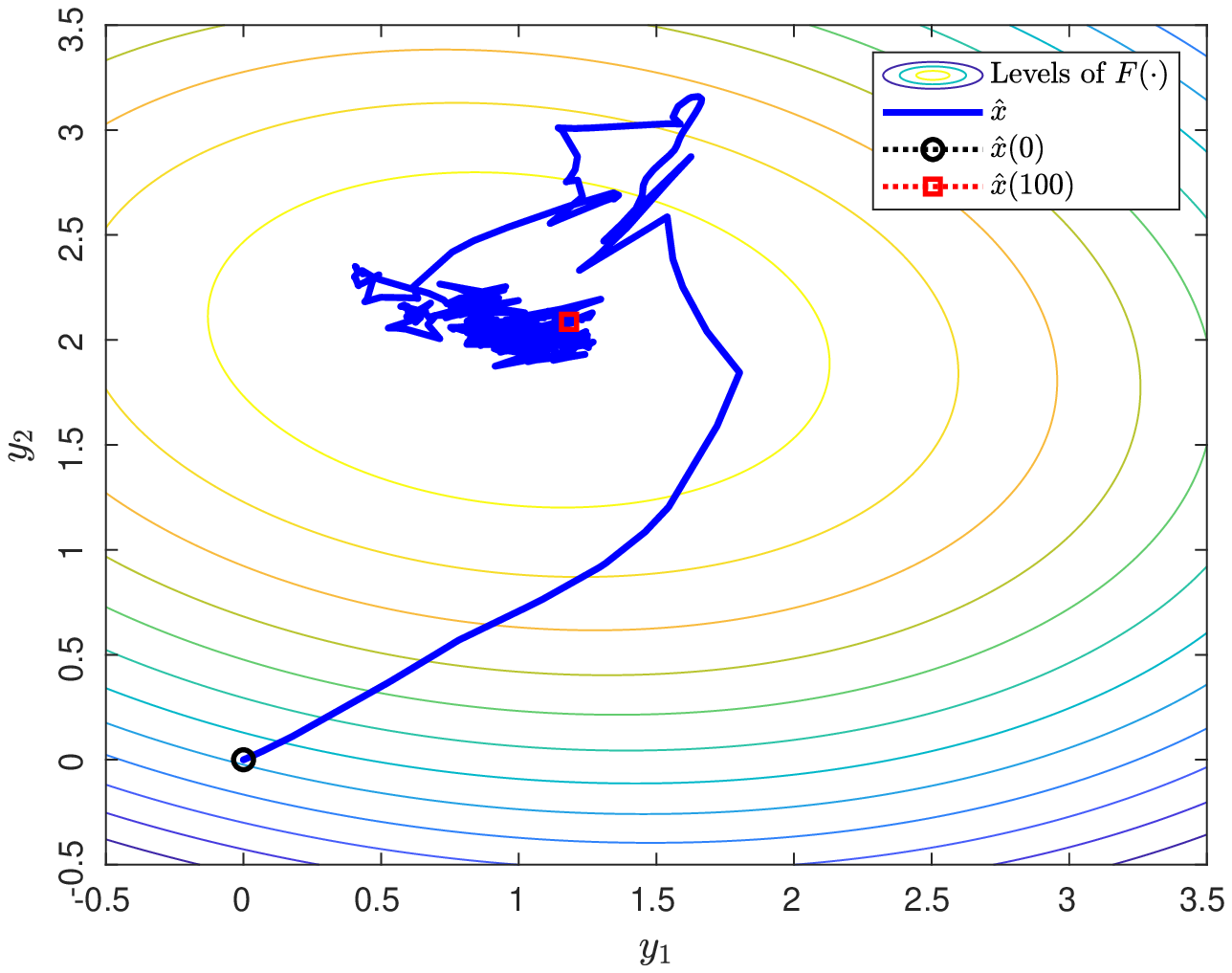}
\end{subfigure}
\caption{ Location estimation for $x(t)=[1 , 2]^T$, $y(t)=[\sin(4t)+\sin(5t), \ \sin(2t)+\sin(3t)]$. Noise in sensing the signal intensity with variance(0.05). The dashed lines and the solid lines represent the actual values and  their estimates, respectively.}
\label{fig:n_nd}
\end{figure*}

Scenario 2: Consider the same conditions in Scenario \ref{fig:nn_nd}, but with white noise with variance(0.05) on $F(t)$ measurement of the sensory agent. Figure \ref{fig:n_nd} displays that the localization is accomplished with some errors scaled with the noise magnitude.

\begin{figure*}
\begin{subfigure}[t]{.33\linewidth}
\includegraphics[width=\linewidth]{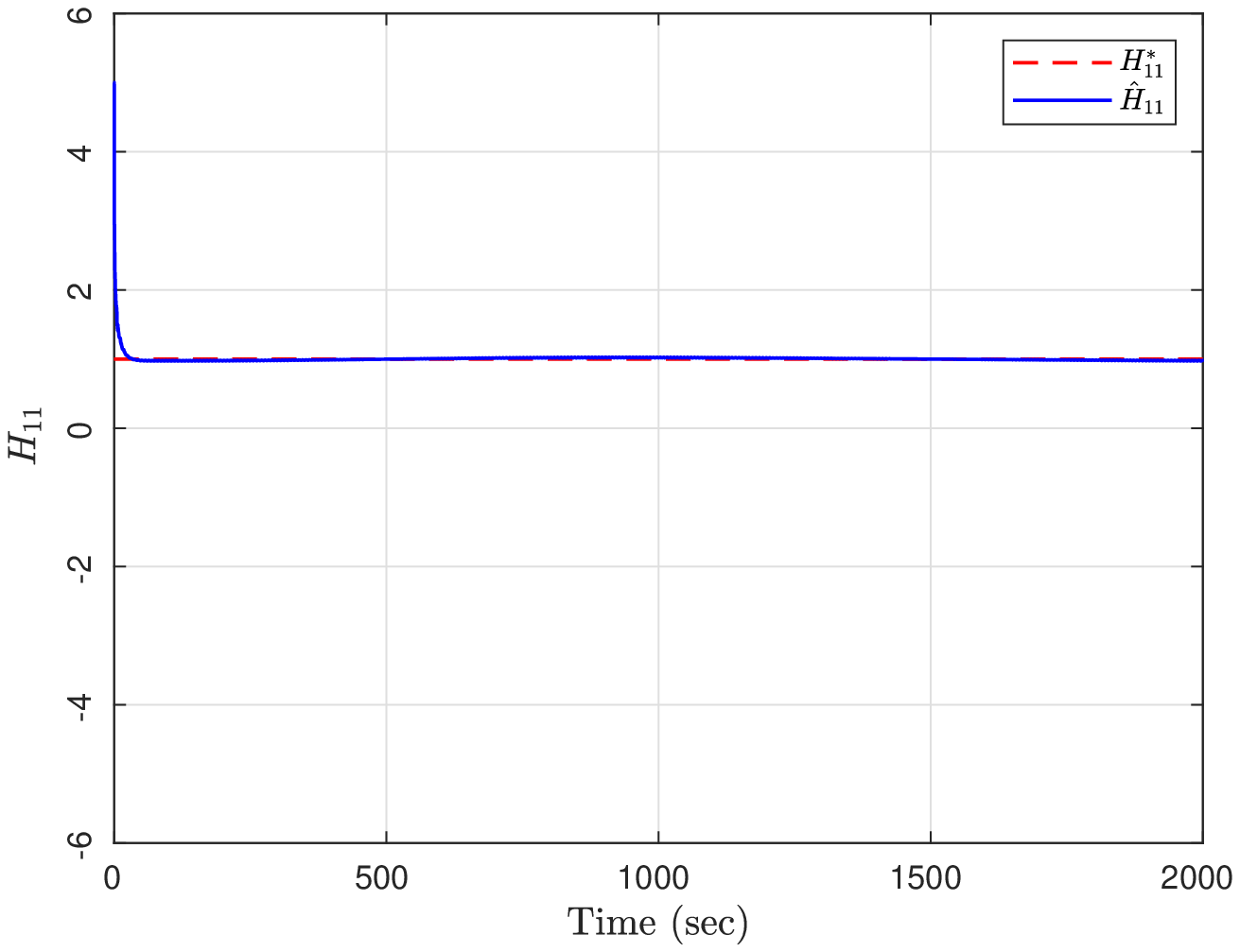}
\end{subfigure}
\begin{subfigure}[t]{.33\linewidth}
\includegraphics[width=\linewidth]{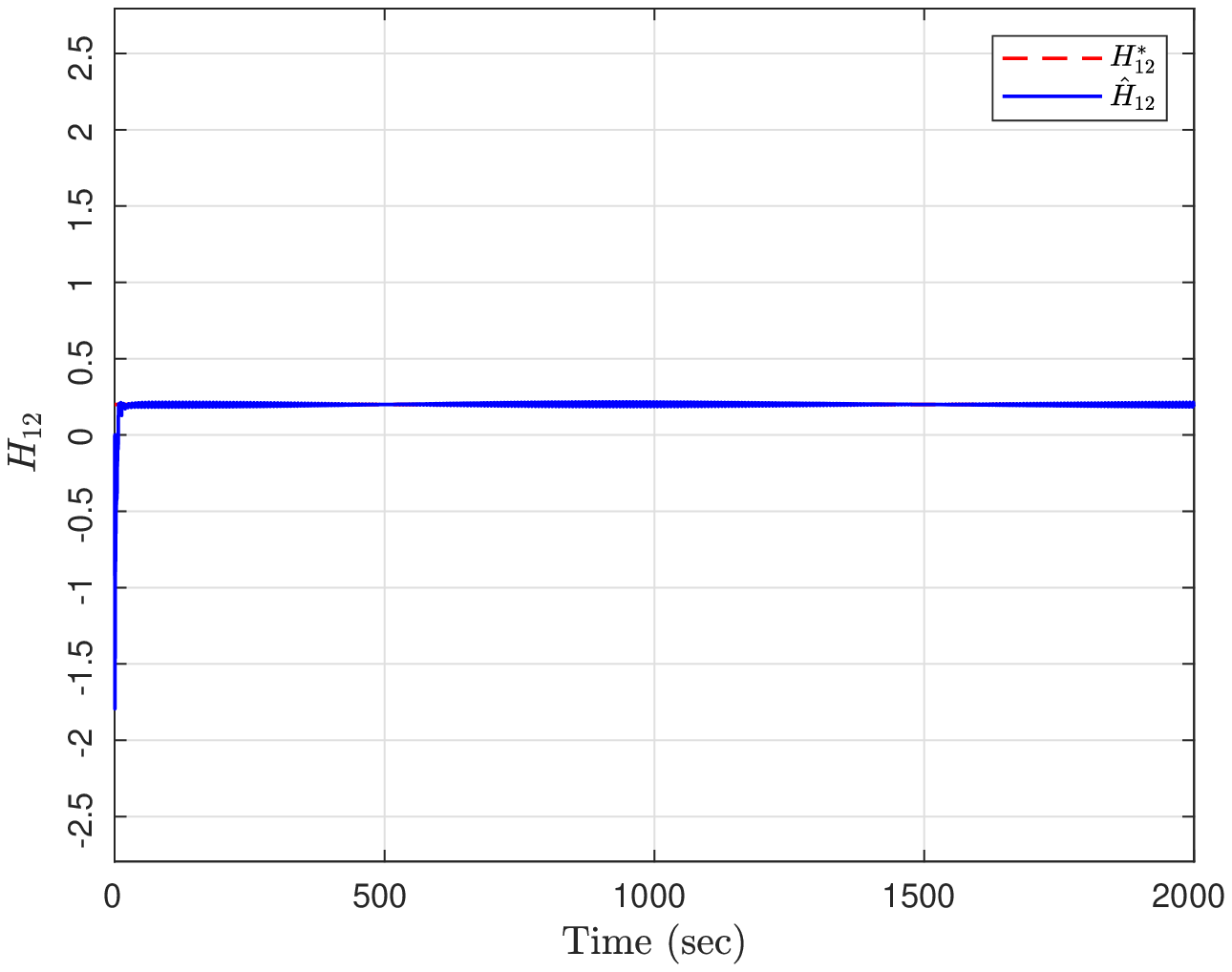}
\end{subfigure}
\begin{subfigure}[t]{.33\linewidth}
\includegraphics[width=\linewidth]{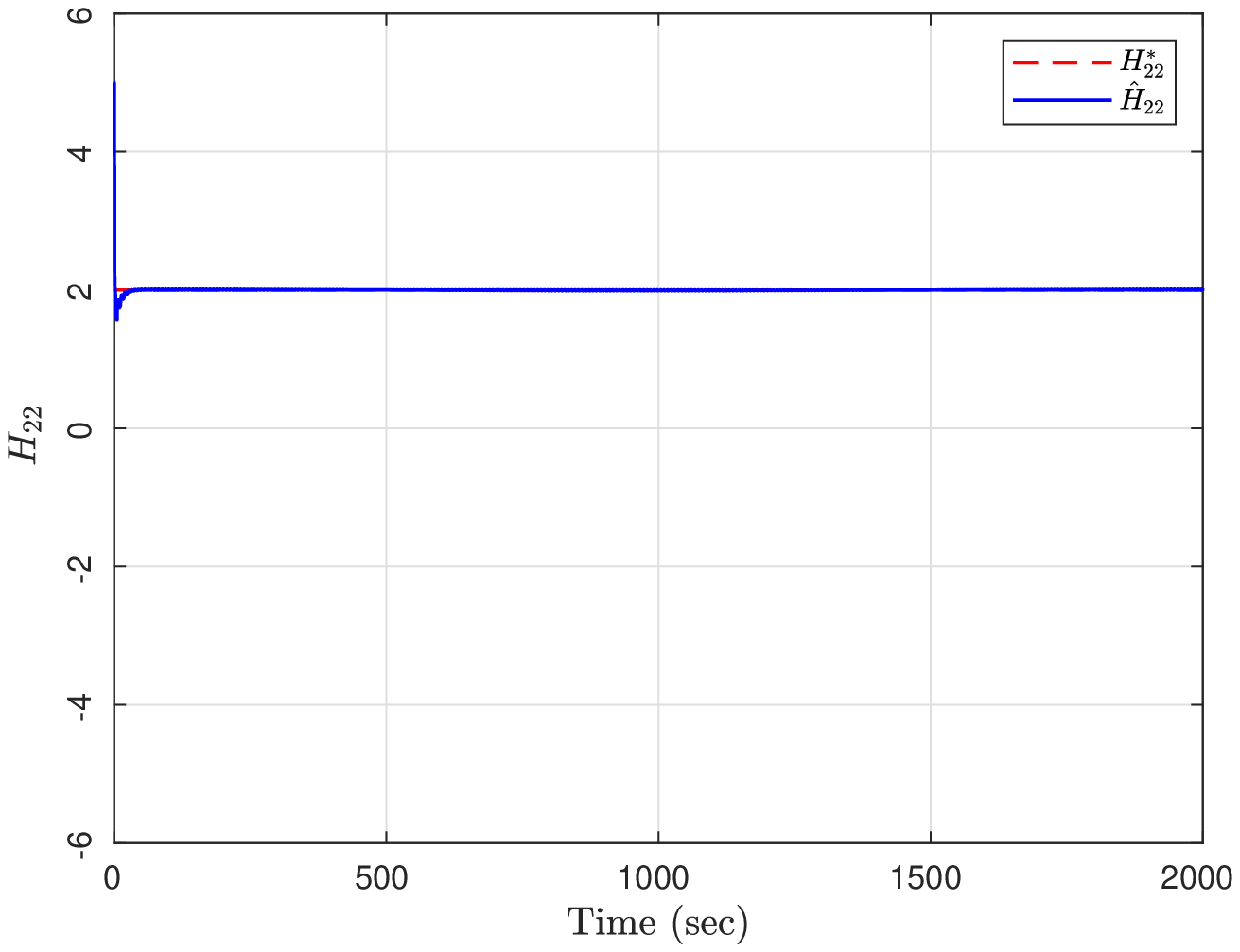}
\end{subfigure}
\begin{subfigure}[t]{.33\linewidth}
\includegraphics[width=\linewidth]{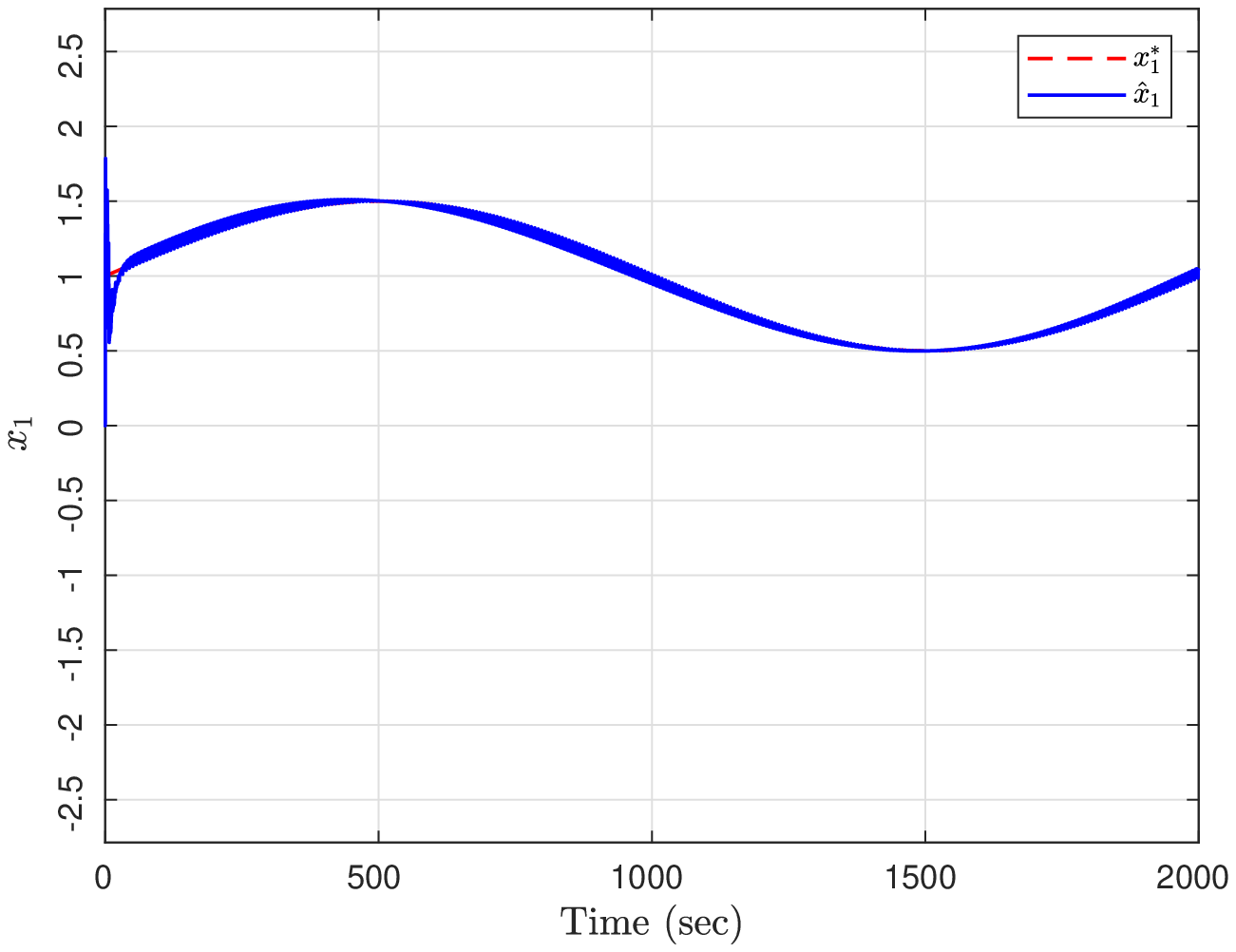}
\end{subfigure}
\begin{subfigure}[t]{.33\linewidth}
\includegraphics[width=\linewidth]{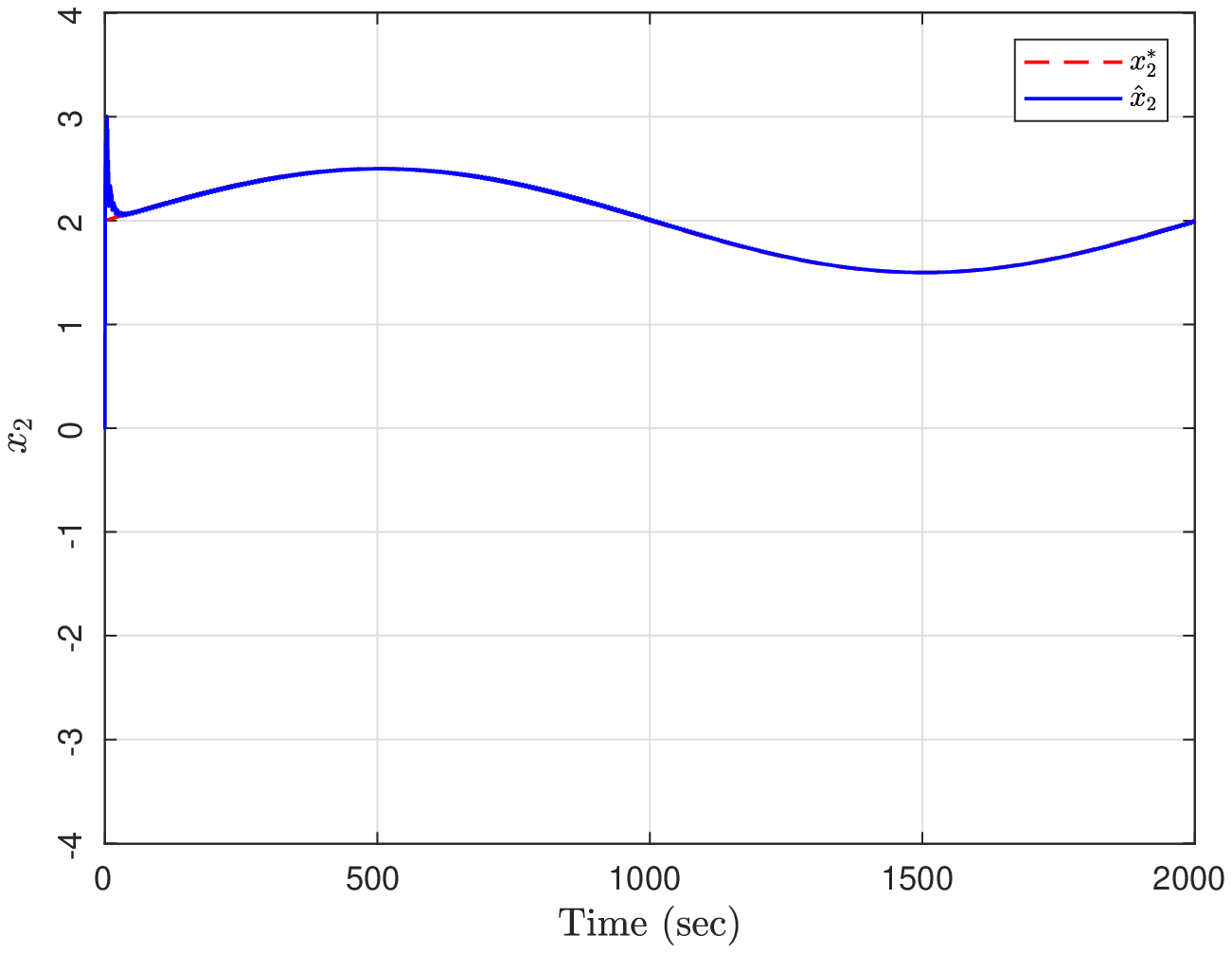}
\end{subfigure}
\begin{subfigure}[t]{.33\linewidth}
\includegraphics[width=\linewidth]{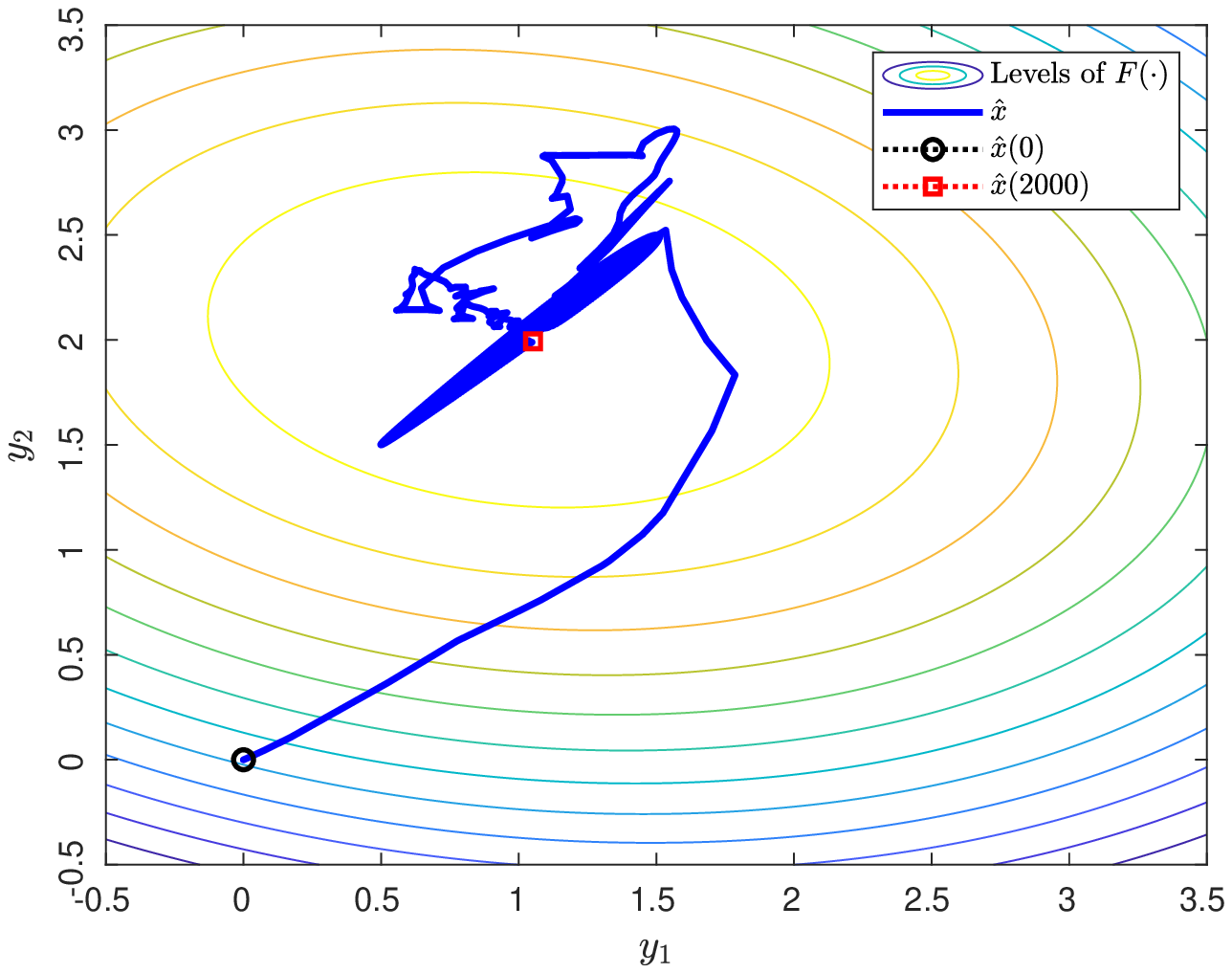}
\end{subfigure}
\caption{ Location estimation for $x(t)=[1+0.5\sin{\frac{\pi}{1000}t}, \ 2+0.5\sin{\frac{\pi}{1000}t}]^T$, $y(t)=[\sin(4t)+\sin(5t), \ \sin(2t)+\sin(3t)]^T$, $a=0.5$. The dashed lines and the solid lines represent the actual values and  their estimates, respectively.}
\label{fig:nn_d}
\end{figure*}

Scenario 3: There is a slow drift movement in the location of extremum as $x(t)=[1+0.5\sin{\frac{\pi}{1000}t}, \ 2+0.5\sin{\frac{\pi}{1000}t}]^T$. As expected from Subsection  \ref{subsec:driftin}, the simulation results in Figure \ref{fig:nn_d} show that the adaptive estimation algorithm in  \eqref{eq:8proj} is applicable for the drift case.

\begin{figure*}
\begin{subfigure}[t]{.33\linewidth}
\includegraphics[width=\linewidth]{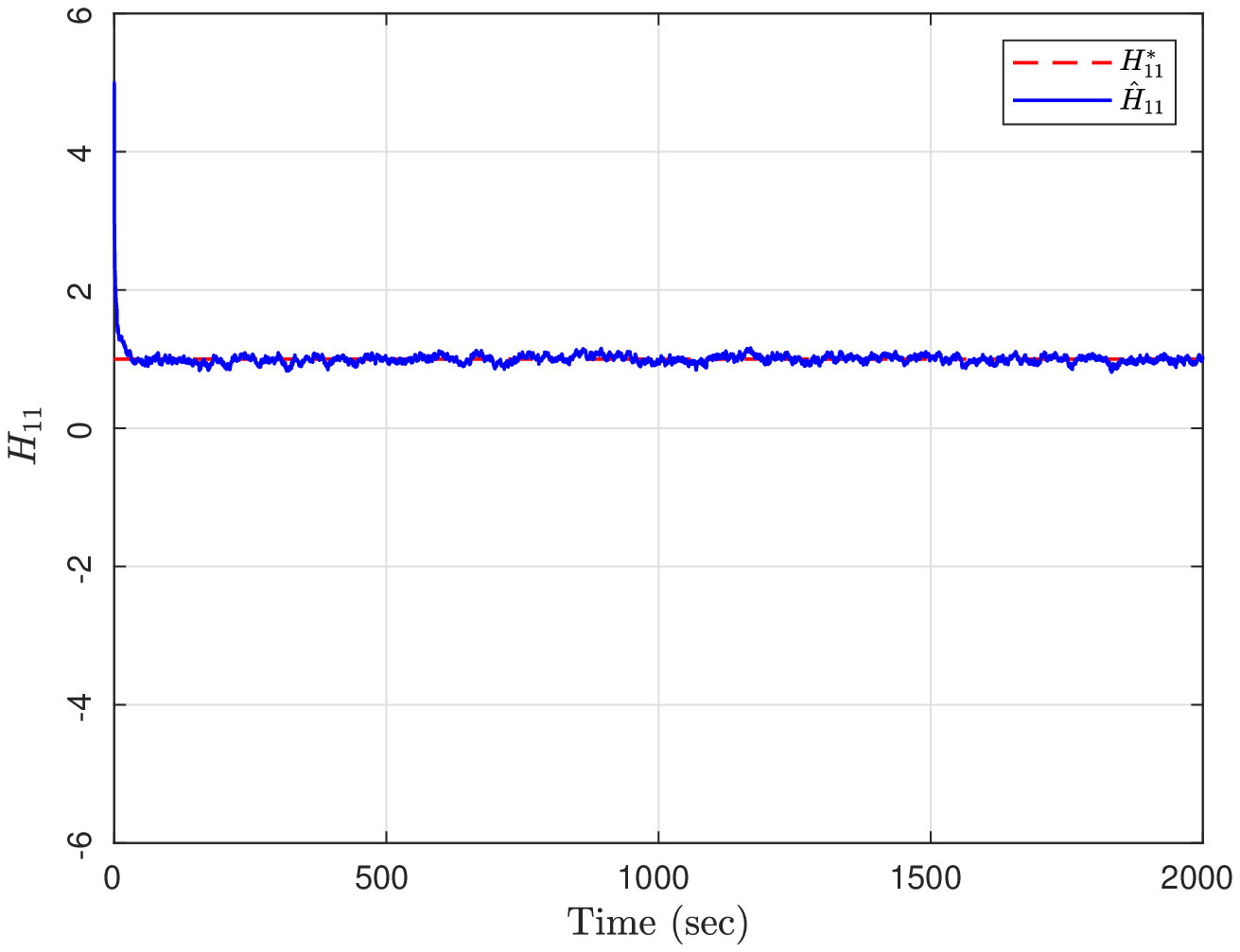}
\end{subfigure}
\begin{subfigure}[t]{.33\linewidth}
\includegraphics[width=\linewidth]{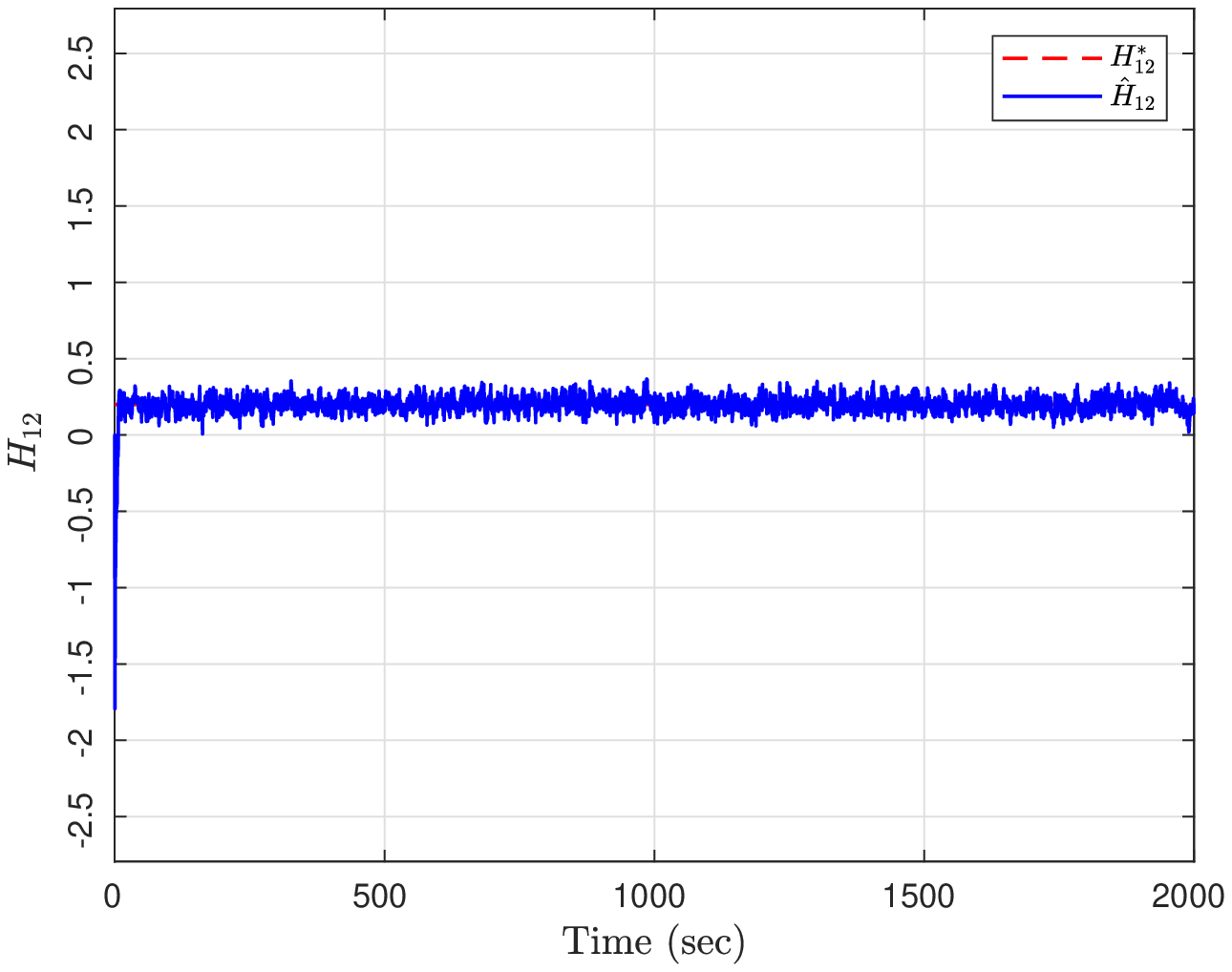}
\end{subfigure}
\begin{subfigure}[t]{.33\linewidth}
\includegraphics[width=\linewidth]{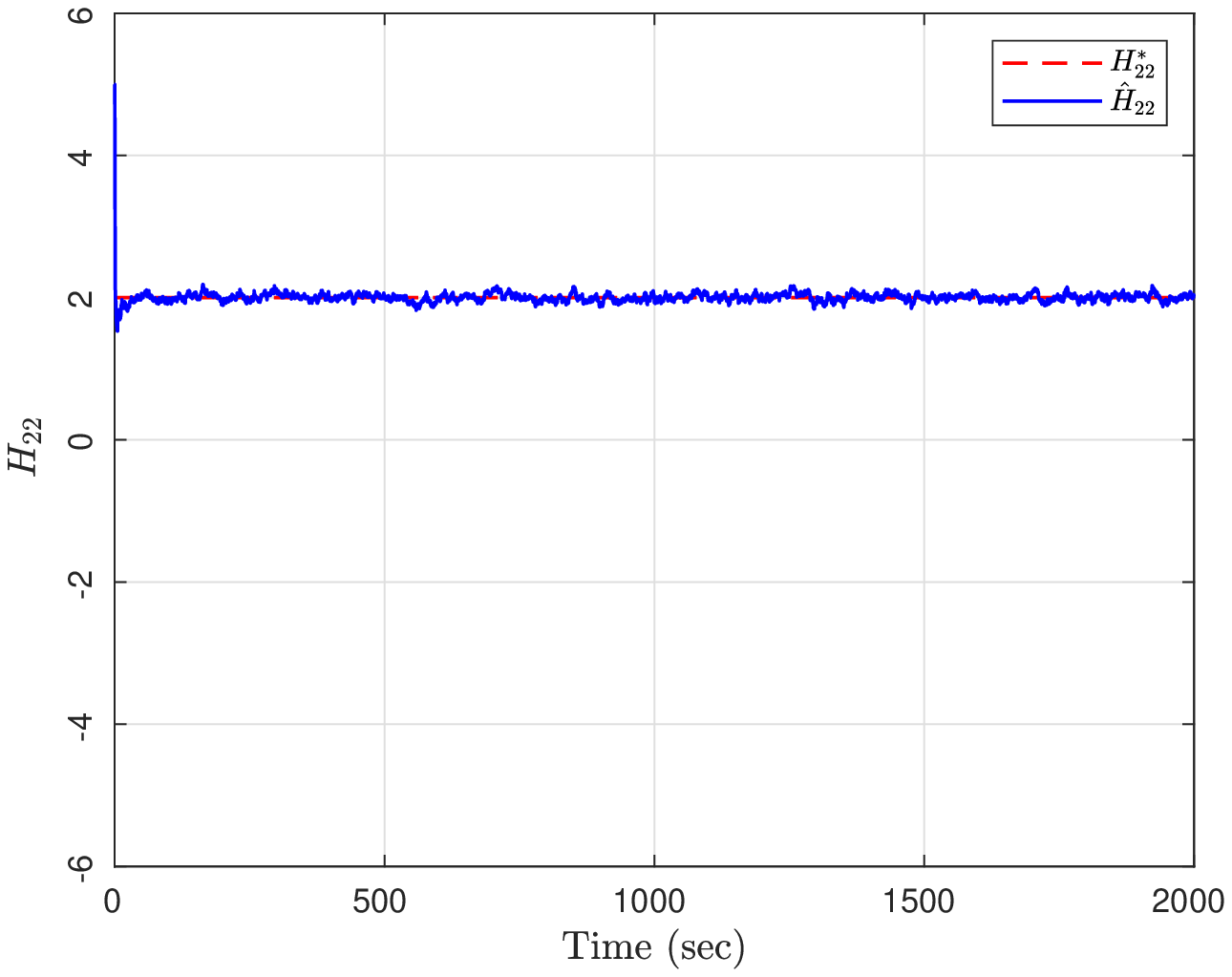}
\end{subfigure}
\begin{subfigure}[t]{.33\linewidth}
\includegraphics[width=\linewidth]{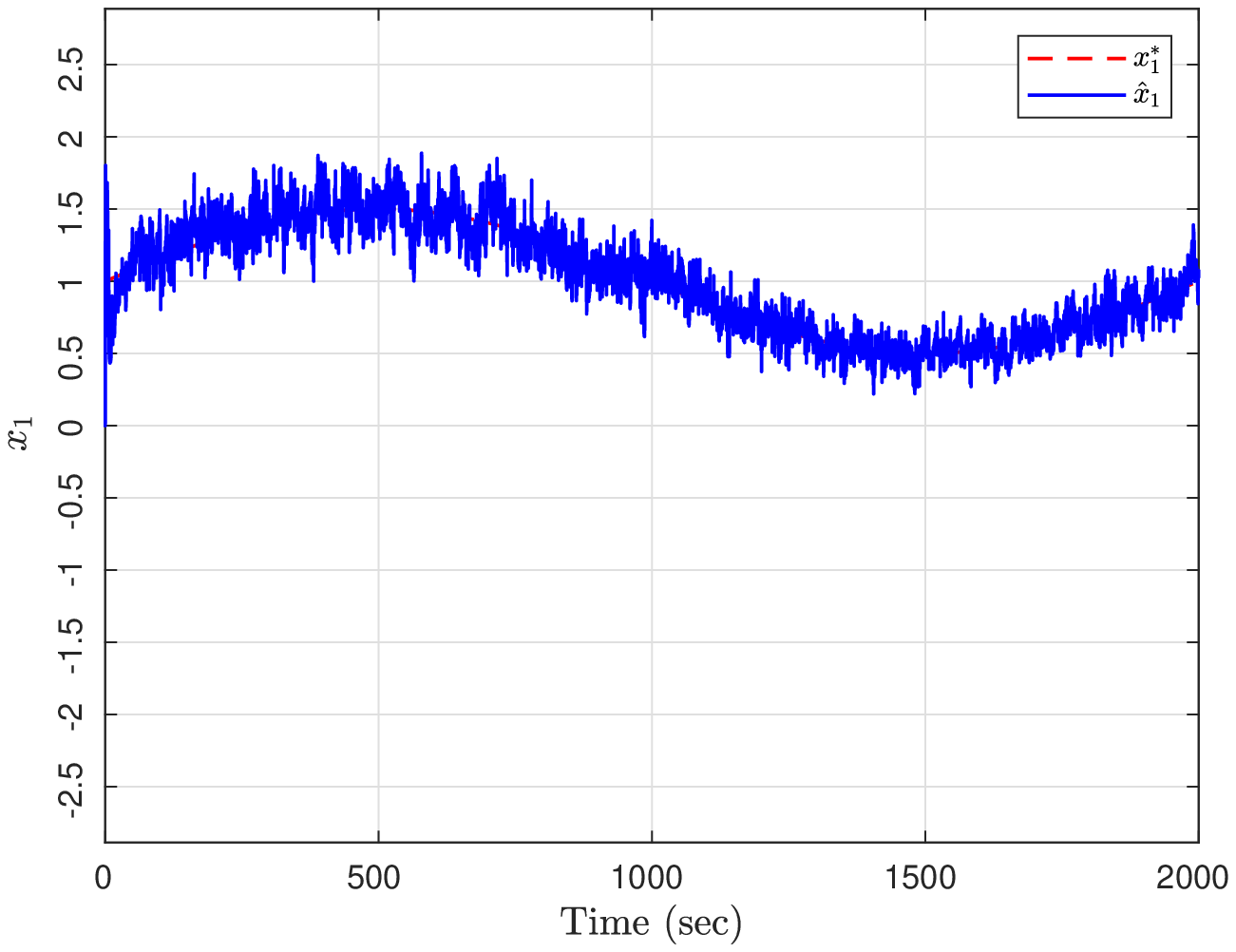}
\end{subfigure}
\begin{subfigure}[t]{.33\linewidth}
\includegraphics[width=\linewidth]{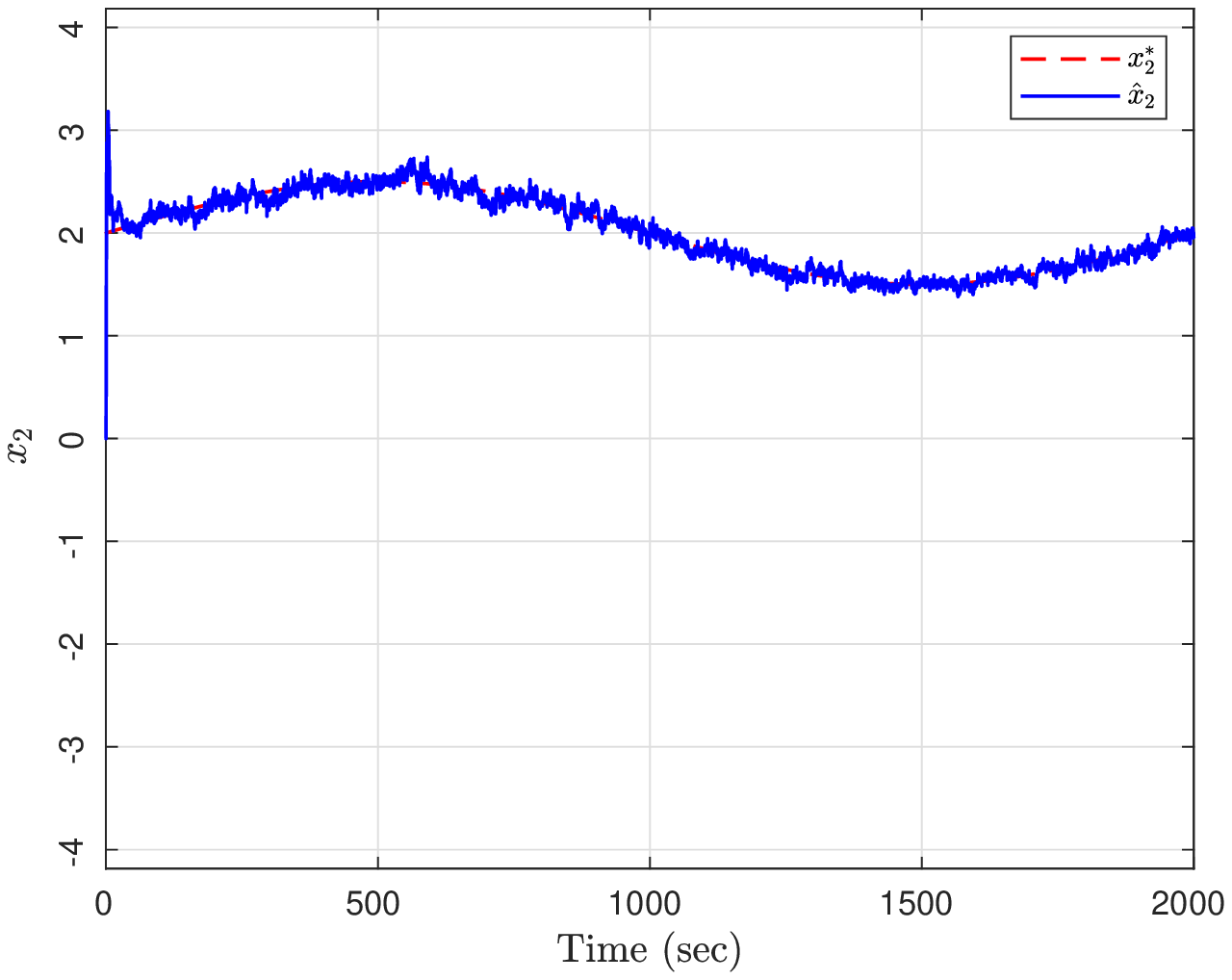}
\end{subfigure}
\begin{subfigure}[t]{.33\linewidth}
\includegraphics[width=\linewidth]{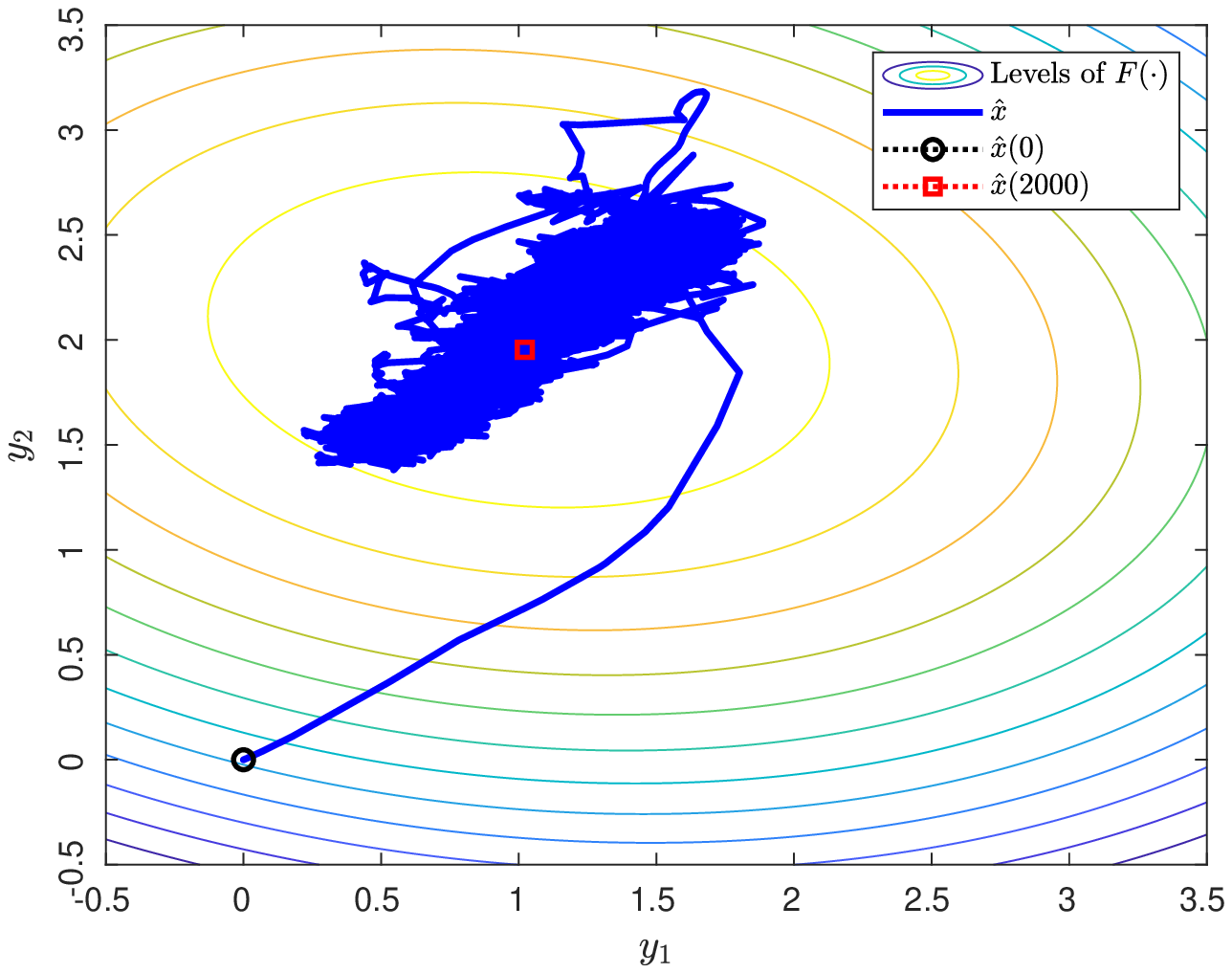}
\end{subfigure}
\caption{ Location estimation for $x(t)=[1+0.5\sin{\frac{\pi}{1000} t}, \ 2+0.5\sin{\frac{\pi}{1000}t}]^T$, $y(t)=[\sin(4t)+\sin(5t), \ \sin(2t)+\sin(3t)]^T$, $a=0.5$. Noise in sensing the signal intensity with variance(0.05). The dashed lines and the solid lines represent the actual values and  their estimates, respectively.}
\label{fig:n_d}
\end{figure*}

Scenario 4: Combine the two circumstances in Scenarios \ref{fig:n_nd} and \ref{fig:nn_d}. There is $F(t)$ measurement noise with variance(0.05) and drift in the location of extremum point as $x(t)=[1+0.5\sin{\frac{\pi}{1000}t}, \ 2+0.5\sin{\frac{\pi}{1000}t}]^T$.  The simulation results in Figure \ref{fig:n_d} demonstrate the adaptive estimation algorithm in \eqref{eq:8proj} works well despite the extremum location drift and noise in sensing.

\section{Conclusion}
\label{sec:Conclusion}
In this paper we have designed an adaptive scheme for Hessian estimation and extremum localization of quadratic signal field functions by a sensory agent measuring the signal intensity. The proposed scheme is effective in extracting more detailed information about such signal fields and utilizing this information in more accurate and faster localization of the extremum. The stability of the proposed adaptive estimation and localization scheme has been proven for both stationary and slowly drifting extremum cases. Simulation results are presented in the presence of realistic measurement noise and drift in extremum location that exhibit the performance of the proposed scheme.

Ongoing and future related research directions include implementing the proposed scheme on autonomous vehicle and cooperative extensions of the design where more than one sensory agent are utilized.

\addtolength{\textheight}{0cm}   % This command serves to balance the column lengths
                                  % on the last page of the document manually. It shortens
                                  % the textheight of the last page by a suitable amount.
                                  % This command does not take effect until the next page
                                  % so it should come on the page before the last. Make
                                  % sure that you do not shorten the textheight too much.

%%%%%%%%%%%%%%%%%%%%%%%%%%%%%%%%%%%%%%%%%%%%%%%%%%%%%%%%%%%%%%%%%%%%%%%%%%%%%%%%

%%%%%%%%%%%%%%%%%%%%%%%%%%%%%%%%%%%%%%%%%%%%%%%%%%%%%%%%%%%%%%%%%%%%%%%%%%%%%%%%

%%%%%%%%%%%%%%%%%%%%%%%%%%%%%%%%%%%%%%%%%%%%%%%%%%%%%%%%%%%%%%%%%%%%%%%%%%%%%%%%

%%%%%%%%%%%%%%%%%%%%%%%%%%%%%%%%%%%%%%%%%%%%%%%%%%%%%%%%%%%%%%%%%%%%%%%%%%%%%%%%

\bibliographystyle{IEEEtran}
\bibliography{adaptiveHess}

% Generated by IEEEtran.bst, version: 1.14 (2015/08/26)
\begin{thebibliography}{10}
\providecommand{\url}[1]{#1}
\csname url@samestyle\endcsname
\providecommand{\newblock}{\relax}
\providecommand{\bibinfo}[2]{#2}
\providecommand{\BIBentrySTDinterwordspacing}{\spaceskip=0pt\relax}
\providecommand{\BIBentryALTinterwordstretchfactor}{4}
\providecommand{\BIBentryALTinterwordspacing}{\spaceskip=\fontdimen2\font plus
\BIBentryALTinterwordstretchfactor\fontdimen3\font minus
  \fontdimen4\font\relax}
\providecommand{\BIBforeignlanguage}[2]{{%
\expandafter\ifx\csname l@#1\endcsname\relax
\typeout{** WARNING: IEEEtran.bst: No hyphenation pattern has been}%
\typeout{** loaded for the language `#1'. Using the pattern for}%
\typeout{** the default language instead.}%
\else
\language=\csname l@#1\endcsname
\fi
#2}}
\providecommand{\BIBdecl}{\relax}
\BIBdecl

\bibitem{mao2009localization}
G.~Mao and B.~Fidan, \emph{Localization Algorithms and Strategies for Wireless
  Sensor Networks}.\hskip 1em plus 0.5em minus 0.4em\relax IGI Global, 2009.

\bibitem{umay2017adaptive}
I.~Umay and B.~Fidan, ``Adaptive wireless biomedical capsule tracking based on
  magnetic sensing,'' \emph{International Journal of Wireless Information
  Networks}, vol.~24, no.~2, pp. 189--199, 2017.

\bibitem{sayed2005network}
A.~H. Sayed, A.~Tarighat, and N.~Khajehnouri, ``Network-based wireless
  location: challenges faced in developing techniques for accurate wireless
  location information,'' \emph{IEEE Signal Processing Magazine}, vol.~22,
  no.~4, pp. 24--40, 2005.

\bibitem{bishop2010optimality}
A.~N. Bishop, B.~Fidan, B.~D. Anderson, K.~Dogancay, and P.~N. Pathirana,
  ``Optimality analysis of sensor-target localization geometries,''
  \emph{Automatica}, vol.~46, no.~3, pp. 479--492, 2010.

\bibitem{niculescu2004positioning}
D.~Niculescu, ``Positioning in ad hoc sensor networks,'' \emph{IEEE Network},
  vol.~18, no.~4, pp. 24--29, 2004.

\bibitem{klukas1998line}
R.~Klukas and M.~Fattouche, ``Line-of-sight angle of arrival estimation in the
  outdoor multipath environment,'' \emph{IEEE Trans. Vehicular Technology},
  vol.~47, no.~1, pp. 342--351, 1998.

\bibitem{cong2002hybrid}
L.~Cong and W.~Zhuang, ``Hybrid tdoa/aoa mobile user location for wideband cdma
  cellular systems,'' \emph{IEEE Trans. Wireless Communications}, vol.~1,
  no.~3, pp. 439--447, 2002.

\bibitem{134479}
W.~A. Gardner and C.-K. Chen, ``Signal-selective time-difference-of-arrival
  estimation for passive location of man-made signal sources in highly
  corruptive environments. i. theory and method,'' \emph{IEEE Trans. Signal
  Processing}, vol.~40, no.~5, pp. 1168--1184, May 1992.

\bibitem{cho2010mobile}
H.~Cho and S.~W. Kim, ``Mobile robot localization using biased
  chirp-spread-spectrum ranging,'' \emph{IEEE Trans. Industrial Electronics},
  vol.~57, no.~8, pp. 2826--2835, 2010.

\bibitem{larsson1996mobile}
U.~Larsson, J.~Forsberg, and A.~Wernersson, ``Mobile robot localization:
  integrating measurements from a time-of-flight laser,'' \emph{IEEE Trans.
  Industrial Electronics}, vol.~43, no.~3, pp. 422--431, 1996.

\bibitem{li2006rss}
X.~Li, ``Rss-based location estimation with unknown pathloss model,''
  \emph{IEEE Trans. Wireless Communications}, vol.~5, no.~12, 2006.

\bibitem{li2002detection}
D.~Li, K.~D. Wong, Y.~H. Hu, and A.~M. Sayeed, ``Detection, classification, and
  tracking of targets,'' \emph{IEEE Signal Processing Magazine}, vol.~19,
  no.~2, pp. 17--29, 2002.

\bibitem{fidan}
S.~Dandach, B.~Fidan, S.~Dasgupta, and B.~Anderson, ``A continuous time linear
  adaptive source localization algorithm, robust to persistent drift,''
  \emph{Systems \& Control Letters}, vol.~58, no.~1, pp. 7--16, 2009.

\bibitem{fidan2013adaptive}
B.~Fidan, S.~Dasgupta, and B.~D. Anderson, ``Adaptive range-measurement-based
  target pursuit,'' \emph{International Journal of Adaptive Control and Signal
  Processing}, vol.~27, no. 1-2, pp. 66--81, 2013.

\bibitem{fidan2015least}
B.~Fidan, A.~Camlica, and S.~Guler, ``Least-squares-based adaptive target
  localization by mobile distance measurement sensors,'' \emph{International
  Journal of Adaptive Control and Signal Processing}, vol.~29, no.~2, pp.
  259--271, 2015.

\bibitem{fidan2015adaptive}
B.~Fidan and I.~Umay, ``Adaptive environmental source localization and tracking
  with unknown permittivity and path loss coefficients,'' \emph{Sensors},
  vol.~15, no.~12, pp. 31\,125--31\,141, 2015.

\bibitem{skobeleva2018planar1}
A.~Skobeleva, B.~Fidan, V.~Ugrinovskii, and I.~R. Petersen, ``Planar
  cooperative extremum seeking with guaranteed convergence using a three-robot
  formation,'' in \emph{Proc. IEEE Conference on Decision and
  Control$\textnormal{, 2018}$, to appear, preprint arXiv:1809.03674}.

\bibitem{brinon2016distributed}
L.~Brinon-Arranz, L.~Schenato, and A.~Seuret, ``Distributed source seeking via
  a circular formation of agents under communication constraints,'' \emph{IEEE
  Trans. Control of Network Systems}, vol.~3, no.~2, pp. 104--115, 2016.

\bibitem{moore2010source}
B.~J. Moore and C.~Canudas-de Wit, ``Source seeking via collaborative
  measurements by a circular formation of agents,'' in \emph{Proc. IEEE
  American Control Conference}, 2010, pp. 6417--6422.

\bibitem{ogren2004cooperative}
P.~Ogren, E.~Fiorelli, and N.~E. Leonard, ``Cooperative control of mobile
  sensor networks: Adaptive gradient climbing in a distributed environment,''
  \emph{IEEE Trans. Automatic Control}, vol.~49, no.~8, pp. 1292--1302, 2004.

\bibitem{krstice1}
M.~Krstic and H.~Wang, ``Stability of extremum seeking feedback for general
  nonlinear dynamic systems,'' \emph{Automatica}, vol.~36, no.~4, pp. 595--601,
  2000.

\bibitem{krstice2}
A.~Ghaffari, M.~Krstic, and D.~Nesic, ``Multivariable newton-based extremum
  seeking,'' \emph{Automatica}, vol.~48, no.~8, pp. 1759--1767, 2012.

\bibitem{krstice3}
S.~Liu and M.~Krstic, ``Newton-based stochastic extremum seeking,''
  \emph{Automatica}, vol.~50, no.~3, pp. 952--961, 2014.

\bibitem{lang2012calculus}
S.~Lang, \emph{Calculus of Several Variables}.\hskip 1em plus 0.5em minus
  0.4em\relax Springer, 2012.

\bibitem{khalil}
H.~Khalil, \emph{Nonlinear Systems}, 3rd~ed.\hskip 1em plus 0.5em minus
  0.4em\relax Prentice Hall, 1996.

\bibitem{ion}
P.~Ioannou{,} and J.~Sun, \emph{Robust Adaptive Control}.\hskip 1em plus 0.5em
  minus 0.4em\relax Prentice Hall, 1996.

\bibitem{ioannou2006adaptive}
P.~Ioannou and B.~Fidan, \emph{Adaptive Control Tutorial}.\hskip 1em plus 0.5em
  minus 0.4em\relax SIAM, 2006, vol.~11.

\bibitem{horn1990matrix}
R.~A. Horn and C.~R. Johnson, \emph{Matrix Analysis}.\hskip 1em plus 0.5em
  minus 0.4em\relax Cambridge University Press, 1990.

\bibitem{anderson1977exponential}
B.~Anderson, ``Exponential stability of linear equations arising in adaptive
  identification,'' \emph{IEEE Trans. Automatic Control}, vol.~22, no.~1, pp.
  83--88, 1977.

\end{thebibliography}

\end{document}